\documentclass[11pt]{article}
\usepackage{mathrsfs}
\usepackage{amsmath}
\usepackage{amsfonts}
\usepackage{amssymb, amsmath, cite}

\def\pa{\partial}
\def\pt{\partial_t}

\newcommand{\N}{\mathbb{N}}
\newcommand{\R}{\mathbb{R}}
\newcommand{\T}{\mathbb{T}}
\newcommand{\Z}{\mathbb{Z}}

\newcommand{\tri}{\|\hspace*{-1pt}|}

\newtheorem{theorem}{Theorem}
\newtheorem{lemma}{Lemma}
\newtheorem{proposition}{Proposition}
\newtheorem{remark}{Remark}

\begin{document}
\begin{center}
{\bf\Large  Stability of non-constant equilibrium solutions for
two-fluid non-isentropic Euler-Maxwell systems arising in plasmas
\\[2mm]}
{Yue-Hong Feng, Xin Li and Shu Wang }  \\[2mm]
{ \it\small  College of Applied Science, Beijing University of
Technology, Ping Le Yuan 100 Beijing 100124, China}\\
 { \it\small Email~: fyh@bjut.edu.cn,\hspace{1mm} lixin91600@163.com;\hspace{1mm}  wangshu@bjut.edu.cn
  }\\
\end{center}
\setlength{\baselineskip}{17pt}{\setlength\arraycolsep{2pt}

\begin{quote}{
{{\bfseries Abstract.} We consider the periodic problem for
 two-fluid  non-isentropic  Euler-Maxwell systems in plasmas. By means of suitable choices of
symmetrizers and an induction argument on the order of the
time-space derivatives of solutions in energy estimates, the global
smooth solution with small amplitude is established near a
non-constant equilibrium solution with asymptotic stability
properties.
This improves the results obtained in \cite{LWF16a} for models with
temperature diffusion terms by using the pressure functions $p^\nu$
in place of  the unknown variables densities $n^\nu$.

 }

}

\end{quote}

\section{Introduction and main results}
Recently  there  are  many  mathematical  researches  on  the
Euler-Maxwell systems  which are used for modeling the motion of
fluid plasmas (see \cite{Ch84,MRS90,RG69} and references theirin).
In this article, we consider the period problem for the two-fluid
non-isentropic compressible Euler-Maxwell system%
\begin{equation}
\label{1.1}
\small\left\{\begin{aligned}
& {{\partial _t}n^\nu + \nabla  \cdot \left( {n^\nu u^\nu} \right) = 0,} \\
&  {\partial _t}({n^\nu }{u^\nu }) + \nabla  \cdot ({n^\nu }{u^\nu }
\otimes {u^\nu }) + \nabla {p^\nu } = {q_\nu }{n^\nu }(E + {u^\nu }
\times B) - {n^\nu }{u^\nu },
 \\
 & {\partial _t} \mathcal {E}^\nu  + \nabla  \cdot\left((\mathcal {E}^\nu +p^\nu
 ) u^\nu \right)
 =  {q_\nu }{n^\nu }{u^\nu } \cdot E - {n^\nu }{{\left| {{u^\nu }}
  \right|}^2} - (\mathcal {E}^\nu - n^\nu e_l),
  \\
& {\partial _t}E - \nabla  \times B =
        {n^e}{u^e}
   -{n^i}{u^i}, \quad \nabla  \cdot E =  {n^i} - {n^e} + b(x), \\
&{\partial _t}B + \nabla  \times E = 0, \quad \nabla
    \cdot B = 0, \quad \nu=e,i, \quad (t,x)\in \R^+ \times \T,
\end{aligned} \right.
\end{equation}
where $\T = (\R/\Z)^3$ denotes a three-dimensional torus and $q_e
=-1 ~~(q_i=1)$ is the charge of electrons (ions). Here, $(n^\nu,
u^\nu, \theta^\nu, E, B)$ $:$ $ \R^+ \times \T \rightarrow $ $ \R^1
\times \R^3 \times \R^1 \times \R^3 \times \R^3, $ and $\nabla$ is
the usual gradient. Physically,  $n^\nu>0$ is the density of fluid,
 $u^\nu $ is the velocity, $\theta^\nu>0$ is the absolute temperature,
 $E$ is
the electric field, and  $B $ is the magnetic field. Functions
 $p^\nu= n^\nu \theta^\nu$, $\displaystyle \mathcal {E}_\nu = n^\nu  e^\nu + \frac{1}{2}
n^\nu |u^\nu|^2 $, $e^\nu= c_* \theta^\nu$ and $b(x)>0$ denote
respectively pressure, total energy, internal energy and
 a doping term.
 The constants $c_* >0$, $\theta_l >0$ and $e_l =c_* \theta_l $ stand for the coefficient of
heat conduction, the back ground temperature and the background
internal energy, respectively.

Next we set $c_* = \theta_l = e_l =1$ for the sake of simplicity.
This is not an essential restriction in the study of global
existence of smooth solutions. Then for $n^\nu>0$, system
\eqref{1.1}  is equivalent to
\begin{equation}
\label{1.1*}
\left\{{\begin{split}
  & {\pa_t}{n^\nu} + \nabla \cdot \left({n^\nu }{u^\nu }\right) = 0, \\
   &
{\partial _t}{u^\nu } + ({u^\nu } \cdot \nabla ){u^\nu } +
{\frac{1}{{{n^\nu }}}\nabla {p^\nu }}
 =
{q_\nu }(E + {u^\nu } \times B) - {u^\nu },
     \\
   &
{\partial _t}{\theta ^\nu } + {u^\nu } \cdot \nabla {\theta ^\nu } +
{\theta ^\nu }\nabla  \cdot {u^\nu } = {\frac{1}{2}{{\left| {{u^\nu
}} \right|}^2}}
 - ({\theta ^\nu } - 1),
 \\
   & {\pa_t}E - \nabla  \times B =    {n^e}{u^e}
   - {n^i}{u^i},~~~~  \nabla  \cdot E = {n^i} - {n^e} + b(x),  \\
   & {\pa_t}B + \nabla  \times E = 0,~~~~ \nabla  \cdot B = 0,
\quad \nu=e,i, \quad (t,x)\in \R^+ \times \T,
\end{split}} \right.
\end{equation}
with the initial condition %
\begin{equation}
\label{1.3}
\left( { n^\nu  , u^\nu, \theta^\nu, E,B} \right){|_{t = 0}} =
\left( n ^\nu_0, u^\nu_0,  \theta^\nu_0,  E_0, B_0 \right), \quad x
\in
 \T,  \quad \nu=e,i,\end{equation}
which satisfies the compatibility condition
\begin{equation}
\label{1.4}
\nabla  \cdot {E_0} = {n^i_0} - {n^e_0} + b(x),~~~~\nabla  \cdot
{B_0} = 0, \quad  x \in
 \T.
\end{equation}

Now we consider equilibrium  solution $ (n ^\nu, u^\nu, \theta^\nu,
E, B) = \left( { \bar n ^\nu \left( x \right), 0, \bar \theta ^\nu,
\bar E\left( x \right),\bar B\left( x \right)} \right)$ to be a
steady-state solution of \eqref{1.1*}. Then we get
\begin{equation}
\label{1.7}
\left\{
\begin{split}
& - \frac{1}{{{{\bar n}^e}}}\nabla {{\bar p}^e} = \bar E
 = \frac{1}{{{{\bar n}^i}}}\nabla {{\bar p}^i},  \quad {{\bar p}^\nu } = {\bar n}^\nu {{\bar \theta }^\nu }, \\
&{{\bar \theta }^\nu } = 1, \\
&\nabla  \times \bar E = 0,  \quad \nabla  \cdot \bar E = {{\bar n}^i} - {{\bar n}^e} + b\left( x \right),   \\
&\nabla  \times \bar B = 0, \quad \nabla  \cdot \bar B = 0.
  \quad  x \in
 \T,
\end{split}  \right.
\end{equation}
which implies that $\bar B$ is a constant vector in $\mathbb{R}^3$.
Moreover, the same analysis as that for Proposition 1.1 in
\cite{LWF16a} gives the existence and uniqueness of a
 smooth periodic solution for (\ref{1.7}).
\begin{lemma}
\label{L1.1}
    Suppose $b =b(x)$ to be a smooth periodic function such that
$b \geq \mbox{const.}
>0$ in $\mathbb{T}$. Then the periodic problem  (\ref{1.7})
has a unique smooth solution  $(\bar{n}^\nu, 0, \bar{\theta}^\nu,
\bar E, \bar B)$ with $\bar \theta ^\nu  = 1$ and $\bar n ^\nu \geq
\mbox{const.}
>0$  in $\mathbb{T}$,  $ \nu=e,i.$
\end{lemma}

For $n^\nu, \theta^\nu
> 0$, \eqref{1.1*} is a nonlinear and symmetrizable
hyperbolic-parabolic system  in the sense of Friedrichs. Then
following the result of Kato \cite{Ka75} and the famous work of
Majda \cite{Ma84}, the periodic problem \eqref{1.1*}-\eqref{1.3}
admits a unique local smooth solution when initial values are
smooth.

\begin{proposition}
\label{P1}
    (Local existence of smooth solutions, see \cite{Ka75,Ma84})
Assume integer $s \geq 3$ and \eqref{1.4} holds. Let $\bar B \in
\R^3$ be any given constant vector and $(\bar{n}^\nu, 0,
\bar{\theta}^\nu, \bar E, \bar B)$ be an equilibrium solution of
\eqref{1.1*} in the sense of Lemma \ref{L1.1}. Suppose
$({n}^\nu_0-\bar{n}^\nu, u^\nu_0, \theta^\nu_0 -1,
E_0-\bar{E},B_0-\bar{B}) \in H^{s}(\T)$ with ${n}^\nu_0,
{\theta}^\nu_0 \geq 2 \kappa$ for some given constant $\kappa>0$.
Then there exists $T > 0$ such that periodic problem
(\ref{1.1*})-(\ref{1.3}) admits a unique smooth solution which
satisfies $n^\nu, \theta^\nu \geq \kappa$ in $[0,T] \times \T$ and
$$ (n^\nu-\bar{n}^\nu, u^\nu, \theta^\nu - 1, E-\bar{E}, B-\bar{B}) \in C^1\big([0,T);H^{s-1}(\T)\big)
\cap C\big([0,T);H^{s}(\T)\big), \quad \nu = e, i.  $$
\end{proposition}

There are mathematical investigations in numerical simulations
\cite{DDS12}, the asymptotic limits with small
parameters\cite{PW08a}, the existence of solutions for Euler-Maxwell
systems. Particularly, some of them are concerned with the global
existence and asymptotic stability of small-amplitude smooth
solutions near constant equilibrium sates. For one-dimensional
simplified isentropic Euler-Maxwell system in which the energy
equation is not contained, the global existence of weak solutions is
proved by the compensated compactness method \cite{CJW00}.
 For the three-dimensional isentropic Euler-Maxwell
systems, the existence of global smooth small solutions to the
Cauchy problem in $\mathbb{R}^3$ is established for $s \ge 3$ and
the asymptotic behaviors of solutions when $s \ge 4$ \cite{UWK12}.
By using suitable choices of symmetrizers and energy estimates,
 the global
existence and the long time behaviors of smooth solutions to the
periodic problem in $\mathbb{T}$ and to the initial value problem in
$\mathbb{R}^3$ for $s \ge 3$ are established \cite{PWG11,Pe12}.
 By high- and low-frequency
decomposition methods, uniform (global) classical solutions to the
initial value problem
 in Chemin-Lerner's spaces with critical regularity is  constructed \cite{Xu11}.
  For $s \ge 4$,
 by the tools of Fourier analysis, the decay rates of global smooth solutions in
$L^q$ with $2 \le q \le \infty$ when the time goes to infinity are
presented \cite{Du11,DLZ12}.
And for $s \ge 6$, the long-time decay rates of global smooth
solutions in $H^{s-2k}(\mathbb{R}^3)$ with $0 \le k \le
\big[s/2\big]$ are also established \cite{UK11}.
For the three-dimensional non-isentropic Euler-Maxwell systems, the
existence of global smooth small solutions to the Cauchy problem in
$\mathbb{R}^3$ is established \cite{FWK11,WFL12}.
  For
the Euler-Maxwell systems without damping, an additional relation
was made to establish such a global existence result for the
one-fluid model \cite{GM11}.
And for the two-fluid case without damping, the global stability of
a constant neutral background is proved \cite{GIP16}, in the sense
that irrotational, smooth and localized perturbations of a constant
background with small amplitude lead to global smooth solutions in
$\R^3$.

All these results above hold when the solutions are close to a
constant equilibrium solution of the Euler-Maxwell systems where
$b(x)$  is a positive constant (for example $b(x) =1$).
When $b(x)$ is a small perturbation of a constant, the Cauchy
problem for compressible Euler-Maxwell systems are considered
\cite{LZ13}, and the time decay rates of smooth solutions are
established.
When $b(x)$ is large, such a stability problem is much more
complicated than before. Recently, motivated by the Guo-Strauss's
work on the study of the damped Euler-Poisson system on the general
bounded domain\cite{GS05}, by employing an induction argument on the
order of the derivatives of solutions, the stabilities of
non-constant equilibrium solutions for the isentropic Euler-Maxwell
systems  \cite{Pe13,FPW15} and non-isentropic Euler-Maxwell systems
with temperature diffusion terms \cite{FWL16a,LWF16a}, respectively.
Very recently, with the help of choosing a new symmetrizer matrix,
the stability of the one-fluid non-isentropic Euler-Maxwell systems
is considered  \cite{LP17}. 
However, there is no result on the stability of non-constant
equilibrium solutions for the two-fluid non-isentropic Euler-Maxwell
systems without temperature diffusion effects so far. The goal here
is to consider this problem.

Now we state the main results of this article.
\begin{theorem}
\label{T1}
  Let $s \geq 3$ and \eqref{1.4} hold. Let $\bar B \in
\R^3$ be any given constant vector and $(\bar{n}^\nu, 0, 1, \bar E,
\bar B)$ be an equilibrium solution of \eqref{1.1*} satisfying
$\bar{n}^\nu \geq \mbox{const.}
>0$ in the sense of Lemma \ref{L1.1}.
Then there exist constants $\delta_0 > 0$ and $C>0$, independent of
any given time $t>0$, such that if
$$  \|(n^\nu_0-\bar{n}^\nu, u^\nu_0,\theta^\nu_0-1, E_0-\bar{E},
 B_0-\bar{B})\|_{s} \leq \delta_0,
\quad \nu=e, i,   $$ where $\|\cdot\|_m$ is the norm of usual
Sobolev spaces $H^m(\T)$, periodic problem (\ref{1.1*})-(\ref{1.3})
has a unique global smooth solution $(n^\nu, u^\nu, \theta^\nu, E,
B)$ satisfying
\begin{equation}
\label{1.12}
\begin{split}
&  \sum\limits_{\nu  = e,i}   ||| \left(  {n^\nu }\left( {t, \cdot }
\right) - {{\bar n}^\nu },{u^\nu }\left( {t, \cdot }
  \right),{\theta ^\nu }\left( {t, \cdot } \right) - 1\right)|||_s^2 + |||\left(E\left( {t, \cdot } \right) - \bar E,B\left( {t, \cdot } \right)
   - \bar B   \right)|||_s^2 \\
  + &  \int_0^t \left( \sum\limits_{\nu  = e,i}  |||\left( {{n^\nu }\left( {\tau , \cdot } \right)
  - {{\bar n}^\nu },{u^\nu }\left( {\tau , \cdot } \right),{\theta ^\nu }\left( {\tau , \cdot }
   \right) - 1} \right)|||_s^2  + |||E\left( {\tau, \cdot } \right) - \bar E|||_{s - 1}^2 \right.\\
& \hspace{1cm} \left.+ |||{\partial _\tau }
  B\left( {\tau , \cdot } \right)|||_{s - 2}^2 + |||{\nabla }B\left( {\tau ,
  \cdot } \right)|||_{s - 2}^2 \right)d\tau  \\
  \le & C\sum\limits_{\nu  = e,i} {} ||\left( {n_0^\nu  - {{\bar n}^\nu },u_0^\nu ,
  \theta _0^\nu  - 1,{E_0} - \bar E,{B_0} - \bar B} \right)||_s^2,
  \quad
  \forall~t
  \ge 0.
\end{split}
\end{equation}
Furthermore,
\begin{equation}
\label{1.13}
\begin{split}
 \mathop{\lim}\limits_{t \to \infty }|||(n^\nu(t)-{\bar n}^\nu, u^\nu(t),
  \theta^\nu(t)-1)|||_{s-1}=0,\quad \nu =e,i,
\end{split}\end{equation}%
\begin{equation}
\label{1.14}
\begin{split}
\mathop {\lim }\limits_{t \to \infty } |||E(t) - \bar{E}|||_{s - 1}
= 0,
\end{split}\end{equation}
and
\begin{equation}
\label{1.15}
\lim_{t \rightarrow +\infty} \big(\tri \pt B(t)\tri_{s-2} + \tri
\nabla B(t)\tri_{s-2}\big) = 0,
\end{equation}%
where $|||\cdot|||_m$ is defined in the next section.
\end{theorem}

\begin{remark}
Obviously, the result above still holds for system \eqref{1.1*} in
case the temperature equation contains temperature diffusion terms
(see \cite{LWF16a}) .
\end{remark}

\begin{remark}
The result in Theorem \ref{T1} for two-fluid non-isentropic
Euler-Maxwell systems still holds for  two-fluid  non-isentropic
Euler-Poisson systems which can be regarded as a special case of the
former systems by setting $B=0$ and $E = - \nabla \Psi$ (see
\cite{Pe13,FPW15}) .
\end{remark}

\begin{remark}
Different from the proof process in \cite{LWF16a}, we choose a new
symmetrizer like \eqref{2.8} here. The effect of temperature
diffusion in non-isentropic Euler-Maxwell equations has been
released successfully by this suitable choices of symmetrizers and
the techniques of using the pressure functions $p^\nu$ in place of
the unknown variables densities $n^\nu$ (see \cite{Ma84,LP17}).
\end{remark}

\begin{remark}
It should be emphasized that the velocity relaxation and temperature
 terms of the considered two-fluid non-isentropic Euler-Maxwell system
here play a key role in the proof of Theorem \ref{T1}. We shall
study in the other forthcoming work the case of non-relaxation for
which the proof is much more complicated to carry out.
\end{remark}

The proof of Theorem \ref{T1} is mainly based on the suitable
choices of symmetrizers and an induction argument on the order of
the time-space derivatives of solutions. These techniques, firstly
employed by Peng \cite{Pe13} in the one-fluid isentropic case and
then extended by Feng-Peng-Wang \cite{FPW15} to the two-fluid
isentropic case, can release the difficulty due to the appearance of
non-constant equilibrium solutions. Besides of these techniques, by
using diffusion effects of temperature,
Feng-Wang-Li\cite{FWL16a,LWF16a} proved the stability of
non-constant equilibrium solution of the periodic problems to
non-isentropic models for $s \ge 6$. Very recently, 
Liu-Peng \cite{LP17} considered the stability of the one-fluid
non-isentropic  models for $s \ge 3$. It should be pointed out that
the techniques of choosing a non-diagonal symmetrizers and making a
change of unknown variables in \cite{LP17} can replace the help of
diffusion effects of temperature used in \cite{FWL16a}.
We remark that the two-fluid non-isentropic Euler-Maxwell systems
are much more complex than both the two-fluid isentropic and the
one-fluid non-isentropic Euler-Maxwell systems because they contain
two different charged fluids energy equations besides the density
and velocity equations. Different from the one-fluid non-isentropic
Euler-Maxwell systems in \cite{LP17}, we shall overcome the
difficulties caused by the coupling of two fluids when we establish
the energy estimates. This can be done by employing new multiplier
functions. Indeed, firstly we introduce a new potential
function $\psi$ such that $\nabla\psi = \bar E-E $, 
 and then another function $\eta^\nu = Q^\nu + q_\nu\psi$ to
establish the estimate of $\partial_t^k Q^\nu$ indirectly (see Lemma
\ref{L3.6}), where $Q^\nu = (\ln p^\nu) - \ln\bar p^\nu$.
 This
yieldes a recurrence relation in Lemma \ref{L3.7} which allows us to
obtain the estimates by induction on $(k,m)$ with $k$ decreasing and
$m$ increasing.
 Then Theorem \ref{T1} follows in the
way by combining these previous estimates with Proposition \ref{P1}
and the standard continuity argument.

We conclude this section by stating the arrangement of the rest of
this article. In Section 2, we recall some useful preliminary Lemmas
and reformulate the periodic problem under consideration. In section
3, detailed energy estimates are established. In section 4, we
complete the proof of Theorems \ref{T1} by using an induction
argument and combining the estimates above.



\section{Preliminaries and Reformulation of periodic problem \eqref{1.1*}-\eqref{1.3}}

\subsection{Preliminaries}
In this subsection, we want to make preliminary works for proving
Theorem \ref{T1}.
 Firstly, let us introduce some notations for the
use throughout this paper. The expression $f \sim g$ means
$\displaystyle \gamma g\leq f\leq \frac{1}{\gamma}g$ for a constant
$0<\gamma<1$. We denote by $\|\cdot\|_s$ the norm of the usual
Sobolev space $H^s(\T)$, and by $\|\cdot\|$ and $\|\cdot\|_{L^p}$
the norms of $L^2(\mathbb T)$ and $L^p(\mathbb T)$, respectively,
where $2 <p \le +\infty$.   We also denote
 by $\langle\cdot,\cdot\rangle$ the inner product over
$L^2({\mathbb T})$.  For a multi-index $\alpha =(
\alpha_1,\alpha_2,\alpha_3)\in\mathbb{N}^3$, we denote
$$\pa^\alpha = \pa_{x_1}^{\alpha_1}
\pa_{x_2}^{\alpha_2}\pa_{x_3}^{\alpha_3}= \pa_{ 1}^{\alpha_1} \pa_{
2}^{\alpha_2}\pa_{ 3}^{\alpha_3},~~~~ |\alpha|=
\alpha_1+\alpha_2+\alpha_3.$$ For
$\alpha=(\alpha_1,\alpha_2,\alpha_3)$ and $\beta =
(\beta_1,\beta_2,\beta_3) \in \N^3$, $\beta \leq \alpha$ stands for
$\beta_j \leq \alpha_j$ for  $j=1,2,3$, and $\beta < \alpha$ stands
for $\beta \leq \alpha$ and $\beta \neq \alpha$. For $T>0$ and
$m\geq 1$, we
 define the Banach space $$\displaystyle {B_{m,T}}\left( \T  \right) = \mathop
  \cap \limits_{k = 0}^m {C^k}\left( {\left[ {0,T} \right],{H^{m - k}}\left( \T  \right)}
  \right),$$ with the norm $$\displaystyle|||f|||_m=  \sqrt{\sum_{k + |\alpha| \leq m}
\|\pa_t^k \pa^\alpha f\|^2} , \quad \forall f \in B_{m,T}\left( \T
\right). $$ Obviously,  $\|\cdot \|_m \leq |||\cdot |||_m.$

 Next, we introduce  three lemmas  which will be used in
the proof of Theorem \ref{T1}.

\begin{lemma}(Poincar\'e inequality, see \cite{Ev98}.)
\label{L2.1}
Let $1 \le p < \infty $ and $\Omega  \subset {\R^3}$ be a bounded
connected open domain with Lipschitz boundary. Then there exists a
constant $C>0$ depending only on $p$ and $\Omega$ such that
$${\left\| {u - {u_\Omega }} \right\|_{{L^p}\left( \Omega  \right)}}
\le C{\left\| {\nabla u} \right\|_{{L^p}\left( \Omega \right)}},
\quad \forall u \in {W^{1,p}}\left( \Omega  \right),$$ %
where
$${u_\Omega } = \frac{1}{{\left| \Omega  \right|}}\int\limits_\Omega
{u\left( x \right)dx} $$
is the average value of $u$ over $ \Omega$.
\end{lemma}

\vspace{3mm}

\begin{lemma}(see \cite{Pe13})
\label{L2.2}
   Let $s \ge 3$ and $u$, $v \in B_{s,T}(\T)$. It holds
\begin{equation}
\label{2.1}
    ||| uv |||_s \leq C ||| u|||_s \, ||| v|||_s.
\end{equation}
\end{lemma}

\begin{lemma}see \cite{Pe13}
\label{L2.3}
    Let $s \ge 3$ and $v \in B_{s,T}(\T)$ satisfying
$\pa_t v = f(x,v,\pa_x v)$, with $f$ being a smooth function such
that $f(x,0,0)=0$. Then for all $t \in [0,T]$, we have
\begin{equation}
\label{2.2}
  \|\pa_t^k \pa^\alpha v(t,\cdot)\| \leq C \|v(t,\cdot)\|_{s},
\quad \forall \, k+|\alpha| \leq s,
\end{equation}
where the positive constant $C$ may depend continuously on
$\|v\|_{s}$.
\end{lemma}


\subsection{Reformulation of periodic problem
\eqref{1.1*}-\eqref{1.3}}

Firstly, we introduce new unknowns variables. Noticing that $p^\nu =
n^\nu \theta^\nu$, then the summation of the first equation
multiplying  $\theta^\nu$ and the third equation in \eqref{1.1*}
multiplying $n^\nu$ gives
\begin{equation}
\label{2.1*}
{\partial _t}{p^\nu } + {u^\nu } \cdot \nabla {p^\nu } + 2{p^\nu
}\nabla  \cdot {u^\nu } = \frac{{{p^\nu }}}{{2{\theta ^\nu
}}}{\left| {{u^\nu }} \right|^2} - \frac{{{p^\nu }}}{{{\theta ^\nu
}}}(\theta ^\nu -1),\quad \nu =e,i.\end{equation}
Let
\begin{equation}
\label{2.1a}
{q^\nu } = \ln {p^\nu }, \quad {\bar q^\nu } = \ln {\bar p^\nu
},\quad \nu =e,i.\end{equation}
Then for $p^\nu >0$, it follows that
\begin{equation}
\label{2.1}
{\partial _t}{q^\nu } + {u^\nu } \cdot \nabla {q^\nu } + 2\nabla
\cdot {u^\nu } = \frac{1}{{2{\theta ^\nu }}}{\left| {{u^\nu }}
\right|^2} - \frac{1}{{{\theta ^\nu }}}(\theta ^\nu-1), \quad \nu
=e,i.\end{equation}

Suppose $\left( n^\nu, u^\nu, \theta^\nu, E,B \right) $ to be a
local smooth solution to the periodic problem
\eqref{1.1*}-\eqref{1.3}. Now, for $\nu=e,i$,  we introduce the
perturbed variables
\begin{equation}
\label{2.2}
{Q^\nu } = {q^\nu } - {\bar q^\nu }, \quad {\Theta^\nu } =
{\theta^\nu } - 1, \quad  F = E - \bar E, \quad G = B - \bar B,
\quad {V^\nu } = \left( {\begin{array}{*{20}{c}}
   {{Q^\nu }} \hfill  \\
   {{u^\nu }} \hfill  \\
   {{\Theta ^\nu }} \hfill  \\
\end{array}} \right), \quad Z = \left( {\begin{array}{*{20}{c}}
   {{V^e}} \hfill  \\
   {{V^i}} \hfill  \\
   F \hfill  \\
   G \hfill  \\
\end{array}} \right).\end{equation}
Substituting these expressions into \eqref{1.1*}, and taking into
account \eqref{1.7}, we obtain
\begin{equation}
\label{2.3}
\left\{{\begin{split}
&{{\partial _t}{Q^\nu } + {u^\nu } \cdot \nabla {Q^\nu } + 2 \nabla  \cdot {u^\nu } + {u^\nu } \cdot \nabla {{\bar q}^\nu } = \frac{1}{{2{\theta ^\nu }}}{{\left| {{u^\nu }} \right|}^2} - \frac{1}{{{\theta ^\nu }}}{\Theta ^\nu },} \hfill  \\
&   {{\partial _t}{u^\nu } + \left( {{u^\nu } \cdot \nabla } \right){u^\nu } + {\theta ^\nu }\nabla {Q^\nu } + {\Theta ^\nu }\nabla {{\bar q}^\nu } = {q_\nu }\left( {F + {u^\nu } \times \left( {\bar B + G} \right)} \right) - {u^\nu },} \hfill  \\
&   {{\partial _t}{\Theta ^\nu } + {u^\nu } \cdot \nabla {\Theta ^\nu } + {\theta ^\nu }\nabla  \cdot {u^\nu } = \frac{1}{2}{{\left| {{u^\nu }} \right|}^2} - {\Theta ^\nu },} \hfill  \\
&   {{\partial _t}F - \nabla  \times G = {n^e}{u^e} - {n^i}{u^i},\nabla  \cdot F = {N^i} - {N^e},} \hfill  \\
&   {{\partial _t}G + \nabla  \times F = 0,\nabla  \cdot G = 0, }
\quad \nu=e,i, \quad (t,x)\in \R^+ \times \T,
\end{split}} \right.
\end{equation}
with the initial condition~:
\begin{equation}
\label{2.7}
Z{|_{t = 0}} = {Z_0} = \left( {Q_0^e,u_0^e,\Theta
_0^e,Q_0^i,u_0^i,\Theta _0^i,{F_0},{G_0}} \right),   \quad x \in
 \T,
\end{equation}
which satisfies the compatibility condition~:
\begin{equation}
\label{2.7*}
\nabla  \cdot {F_0} = {N^i_0} - {N^e_0}, ~~~~\nabla \cdot {G_0} = 0,
\quad x \in
 \T.
\end{equation}
Here ${N^\nu _0} = {n^\nu _0} - {\bar n^\nu }$, $Q_0^\nu 
= \ln \left( {n_0^\nu \theta _0^\nu } \right) - \ln \left( {{{\bar
n}^\nu } } \right)$, ${\Theta^\nu _0} = {\theta^\nu _0} - 1$, ${F_0}
= {E_0} - \bar E$ and ${G_0} = {B_0} - \bar B$.

A direct computation gives
\begin{equation}
\label{2.4}
 {\begin{split} {N^\nu } = {n^\nu } - {\bar n^\nu } = \frac{{{p^\nu
}}}{{{\theta ^\nu }}} - \frac{{{{\bar p}^\nu }}}{{{{\bar \theta
}^\nu }}} = \frac{{{e^{{q^\nu }}}}}{{{\theta ^\nu }}} -
\frac{{{e^{{{\bar q}^\nu }}}}}{{{{\bar \theta }^\nu }}} \sim Q^\nu +
\Theta ^\nu , \quad \nu=e,i,
\end{split}}
\end{equation}
which implies that $N^\nu$ can be regarded as a function of $Q^\nu$
and $\Theta^\nu$ with order one.

Next, we can also rewrite the non-isentropic Euler equations of
\eqref{2.3} in the matrix form~: %
\begin{equation}
\label{2.5}
{\partial _t}{V^\nu } + \sum\limits_{j = 1}^3 {} A_j^\nu \left(
{{u^\nu },{\theta ^\nu }} \right){\partial _j}{V^\nu } + {L^\nu
}\left( x \right){V^\nu } = {K^\nu }\left( {{u^\nu },{\theta ^\nu
},F,G,x} \right),\quad \nu =e,i,
\end{equation}
Supplemented by the Maxwell equations
\begin{equation}
\label{2.6}
 \left\{\begin{split}
   {{\partial _t}F - \nabla  \times G = {n^e}{u^e} - {n^i}{u^i}, \quad \nabla  \cdot F = {N^i} - {N^e},} \hfill  \\
   {{\partial _t}G + \nabla  \times F = 0, \quad \nabla  \cdot G = 0,\quad (t,x)\in \R^+ \times \T,
   }
\end{split} \right.
\end{equation}
with
\begin{equation}
\label{2.6a}
 A_j^\nu \left( {{u^\nu },{\theta ^\nu }} \right) = \left( {\begin{array}{*{20}{c}}
   \vspace{0.3cm} {u_j^\nu } & {2e_j^T} & 0  \\
  \vspace{0.3cm} {{\theta ^\nu }{e_j}} & {u_j^\nu {I_3}} & 0  \\
   0 & {{\theta ^\nu }e_j^T} & {u_j^\nu }  \\
\end{array}} \right), ~~~~j = 1,2,3,\quad \nu =e,i,
\end{equation}
\begin{equation}
\label{2.6b}
{{ }}{{} {L}^\nu }\left( x \right) = \left( {\begin{array}{*{20}{c}}
  \vspace{0.3cm} 0 & {{{\left( {\nabla {{\bar q}^\nu }} \right)}^T}} & 0  \\
 \vspace{0.3cm}  0 & 0 & {\nabla {{\bar q}^\nu }}  \\
   0 & 0 & 0  \\
\end{array}} \right),\quad \nu =e,i,
\end{equation}
\begin{equation}
\label{2.6c}
{K^\nu }\left( {{u^\nu },{\theta ^\nu },F,G,x} \right) = \left(
{\begin{array}{*{20}{c}}
  \vspace{0.3cm} {\displaystyle\frac{1}{{2{\theta ^\nu }}}{{\left| {{u^\nu }} \right|}^2} - \displaystyle\frac{1}{{{\theta ^\nu }}}{\Theta ^\nu }}  \\
  \vspace{0.3cm} {{q_\nu }\left( {F + {u^\nu } \times \left( {\bar B + G} \right)} \right) - {u^\nu }}  \\
   {\displaystyle\frac{1}{2}{{\left| {{u^\nu }} \right|}^2} - {\Theta ^\nu }}  \\
\end{array}} \right),\quad \nu =e,i,
\end{equation}
where $\left( {{e_1},{e_2},{e_3}} \right)$ denotes the canonical
basis of $\R^3$, ${I_3}$ denotes the $3\times 3$ unit matrix and we
use $[\cdot]^T$  to denote the transpose of a vector.

Obviously, system \eqref{2.5} for $V^\nu$ is symmetrizable
hyperbolic when both $\theta^\nu$ and $p^\nu = n^\nu\theta^\nu $ are
positive. In fact, we can choose a symmetric and positive definite
matrix as
\begin{equation}
\label{2.8}
A_0^\nu ({p^\nu },{\theta ^\nu }) = \left( {\begin{array}{*{20}{c}}
  \vspace{0.3cm} {{p^\nu }} & 0 & { -\displaystyle \frac{{{p^\nu }}}{{{\theta ^\nu }}}}  \\
  \vspace{0.3cm} 0 & {\displaystyle\frac{{{p^\nu }}}{{{\theta ^\nu }}}{I_3}} & 0  \\
   { - \displaystyle\frac{{{p^\nu }}}{{{\theta ^\nu }}}} & 0 & {\displaystyle\frac{{2{p^\nu }}}
   {{{{\left| {{\theta ^\nu }} \right|}^2}}}}  \\
\end{array}} \right),\quad \nu =e,i,
\end{equation}
which implies that
\begin{equation}
\label{2.8a}
\tilde A_j^\nu ({p^\nu },{u^\nu },{\theta ^\nu }) = A_0^\nu ({p^\nu
},{\theta ^\nu }) A_j^\nu \left( {{u^\nu },{\theta ^\nu }} \right) =
\left( {\begin{array}{*{20}{c}}
   \vspace{0.3cm} {{p^\nu }u_j^\nu } & {{p^\nu }e_j^T} & { - \displaystyle\frac{{{p^\nu }}}{{{\theta ^\nu }}}u_j^\nu }  \\
   \vspace{0.3cm} {{p^\nu }{e_j}} & {\displaystyle\frac{{{p^\nu }}}{{{\theta ^\nu }}}u_j^\nu {I_3}} & 0  \\
   { - \displaystyle\frac{{{p^\nu }}}{{{\theta ^\nu }}}u_j^\nu } & 0 & {\displaystyle\frac{{2{p^\nu }}}{{{{
   \left| {{\theta ^\nu }} \right|}^2}}}u_j^\nu }  \\
\end{array}} \right),\quad \nu =e,i,\end{equation}
 is symmetric, and then system \eqref{2.5} is symmetrizable
hyperbolic when $n^\nu>0$ and $\theta^\nu>0$. Moreover, we have%
\begin{equation}
\label{2.8b}
A_0^\nu ({p^\nu },{\theta ^\nu }){L^\nu }(x) = \left(
{\begin{array}{*{20}{c}}
   \vspace{0.3cm} 0 & {{p^\nu }{{\left( {\nabla {{\bar q}^\nu }} \right)}^T}} & 0  \\
  \vspace{0.3cm}  0 & 0 & {\displaystyle\frac{{{p^\nu }}}{{{\theta ^\nu }}}\nabla {{\bar q}^\nu }}  \\
   0 & { - \displaystyle\frac{{{p^\nu }}}{{{\theta ^\nu }}}{{\left( {\nabla {{\bar q}^\nu }} \right)}^T}} & 0  \\
\end{array}} \right),\quad \nu =e,i.\end{equation}
Then it follows from $\nabla {\bar q^\nu } =
\displaystyle\frac{1}{{{{\bar p}^\nu }}}\left( {\nabla {{\bar p}^\nu
}} \right)$ that the following matrix
\begin{equation}
\label{2.8c}
\begin{split}
 {B^\nu }({p^\nu },{u^\nu },{\theta ^\nu },x)
  = & \sum\limits_{j = 1}^3 {} {\partial _j}\tilde A_j^\nu ({p^\nu },{u^\nu },
  {\theta ^\nu }) - 2A_0^\nu ({p^\nu },{\theta ^\nu }){L^\nu }(x) \\
  = & \left( {\begin{array}{*{20}{c}}
  \vspace{0.3cm} {\nabla  \cdot \left( {{p^\nu }{u^\nu }} \right)} & {{{\left( {\nabla {p^\nu }} \right)}^T} - 2
  \displaystyle\frac{{{p^\nu }}}{{{{\bar p}^\nu }}}{{\left( {\nabla {{\bar p}^\nu }} \right)}^T}} & {
  - \nabla  \cdot \left( {\displaystyle\frac{{{p^\nu }{u^\nu }}}{{{\theta ^\nu }}}} \right)}  \\
 \vspace{0.3cm}  {\nabla {p^\nu }} & {\nabla  \cdot \left( {\displaystyle\frac{{{p^\nu }{u^\nu }}}{{{\theta ^\nu }}}}
  \right){I_3}} & { - 2\displaystyle\frac{{{p^\nu }}}{{{\theta ^\nu }{{\bar p}^\nu }}}\left( {\nabla {{\bar p}^\nu }} \right)}  \\
   { - \nabla  \cdot \left( {\displaystyle\frac{{{p^\nu }{u^\nu }}}{{{\theta ^\nu }}}} \right)} & {2
   \displaystyle\frac{{{p^\nu }}}{{{\theta ^\nu }{{\bar p}^\nu }}}{{\left( {\nabla {{\bar p}^\nu }} \right)}^T}}
    & {2\nabla  \cdot \left( {\displaystyle\frac{{{p^\nu }{u^\nu }}}{{{{\left| {{\theta ^\nu }} \right|}^2}}}} \right)}  \\
\end{array}} \right),\quad \nu =e,i, \\
 \end{split}\end{equation}
is antisymmetric at the point $({p^\nu },{u^\nu },{\theta ^\nu }) =
({\bar p^\nu },0,{\bar \theta ^\nu }) = ({\bar p^\nu },0,1).$

From now on, let $T> 0$ and $Z$ be a smooth solution of
\eqref{2.5}-\eqref{2.6} with the initial condition \eqref{2.7}
defined on time interval $[0,T]$. We set
\begin{equation}
\label{3.1a}
\omega_T = \mathop {\sup }\limits_{0 \le t \le T} |||Z\left( t
\right)||{|_{s}}.
\end{equation}
From the continuous embedding $H^{s-1}(\T) \hookrightarrow
L^\infty(\T)$ for $s \geq 3$, there is a constant $C>0$ such that
$$   \|f\|_{L^\infty} \leq C \|f\|_{s-1}, \quad \forall \, f \in H^{s-1}(\T).  $$
In the following, we always suppose that integer $s \ge 3$ and
$\omega_T $ is small enough, as a consequence,
it is easy to see that

\begin{equation}
\label{3.1}
\begin{split}
0 < \mbox{const.}  \le \frac{3}{4}{{\bar n}^\nu } \le {n^\nu }=
\bar{n}^\nu + N^\nu \le \frac{4}{3}{{\bar n}^\nu }, \quad
\frac{3}{4} \le {\theta ^\nu } = 1 + \Theta^\nu \le \frac{4}{3}, \quad \nu =e,i, 
\end{split}
\end{equation}
and then
\begin{equation}
\label{3.1b}
\begin{split}
0 < \mbox{const.} \le \frac{9}{{16}}{{\bar n}^\nu } \le {p^\nu } =
{n^\nu }{\theta ^\nu } \le \frac{{16}}{9}{{\bar n}^\nu },\quad \nu
=e,i,
\end{split}
\end{equation}
\section{Energy estimates}

In this section we begin to use the normal energy method to obtain
some uniform-in-time a priori estimates for smooth solutions to the
periodic problem \eqref{2.5}-\eqref{2.6} with the initial condition
\eqref{2.7}.   In the first subsection, we establish energy
estimates with dissipation estimates of velocity $u^\nu$ and
temperature $\Theta^\nu$. In the second subsection, we establish a
recurrence relationship for proving Theorem \ref{T1} by combining
the dissipation estimates for density $N^\nu$, pressure variables
$Q^\nu$ and electric field $F$.

\subsection{Dissipation energy estimates for velocity $u^\nu$ and temperature $\Theta^\nu$.}

Assume $k \in \N$ and $\alpha \in \N^3$ with $1 \leq k+|\alpha| \leq
s$.  Applying $\pa_t^k \pa^\alpha$ to (\ref{2.5}), we have
\begin{equation}
\label{3.3}
\begin{split}
{\partial _t}V_{k,\alpha }^\nu  + \sum\limits_{j = 1}^3 {A_j^\nu
\left( {{u^\nu },{\theta ^\nu }} \right)} {\partial _j}V_{k,\alpha
}^\nu  + {L^\nu }\left( x \right)V_{k,\alpha }^\nu  =  K_{k,\alpha
}^\nu
 + g_{k,\alpha }^\nu , \quad \nu = e,
i,
 \end{split}\end{equation}
where $$ V^\nu_{k, \alpha} = \pa_t^k{\pa^\alpha }{V^\nu }, \quad
K^\nu_{k, \alpha} = \pa_t^k{\pa^\alpha }{K^\nu },$$ 
and
\begin{equation}
\label{3.5}
\begin{split} g_{k,\alpha }^\nu  =& \sum\limits_{j = 1}^3 {} \left(
{A_j^\nu \left( {{u^\nu },{\theta ^\nu }} \right){\partial
_j}V_{k,\alpha }^\nu  -
\partial _t^k{\partial ^\alpha }\left( {A_j^\nu \left( {{u^\nu
},{\theta ^\nu }} \right){\partial _j}{V^\nu }} \right)} \right)\\
& + {L^\nu }\left( x \right)V_{k,\alpha }^\nu  - \partial
_t^k{\partial ^\alpha }\left( {{L^\nu }\left( x \right)}V^\nu
\right), \quad \nu = e,i.
 \end{split}\end{equation}
For the Maxwell equations, we also have
\begin{equation}
\label{3.4}
\left\{\begin{split}
  & {{\pa_t} F_{k,\alpha} -
    \nabla  \times    G_{k,\alpha}
    = \pa_t^k{\pa^\alpha }\left(  n^e   u^e -  n^i   u^i
    \right), ~~
    \nabla  \cdot F_{k,\alpha} =
     N^i_{k, \alpha} - N^e_{k, \alpha} ,}  \\
   &   \pa_t G_{k,\alpha} +
   \nabla  \times F_{k,\alpha} =
   0,~~
   \nabla  \cdot G_{k,\alpha} = 0,
 \end{split} \right.\end{equation}
with $F_{k,\alpha} = \pa_t^k \pa^\alpha F$, $G_{k,\alpha} = \pa_t^k
\pa^\alpha G$, $N^\nu_{k,\alpha} = \pa_t^k \pa^\alpha N^\nu$,
$\forall~\nu = e,i$, and etc.
\begin{lemma}
\label{L3.3}
Suppose that the conditions of Theorem \ref{T1} hold and $\omega_T$
is small enough independent of $T$, then there exists a positive
constant $C_0$ such that, for all $k \in \N$ and $\alpha \in \N^3$
with $|\alpha| \geq 1$ and $k+|\alpha| \leq s$, we have
\begin{equation}
\label{3.6}
\begin{split}
 & \frac{d}{{dt}}\left( {\sum\limits_{\nu  = e,i} {} \left\langle {A_0^\nu ({p^\nu },{\theta ^\nu
})V_{k,\alpha }^\nu ,V_{k,\alpha }^\nu } \right\rangle  + ||{F_{k,\alpha }}|{|^2} + ||{G_{k,\alpha }}|{|^2}} \right) + {C_0}\sum\limits_{\nu  = e,i} {} \left( {||u_{k,\alpha }^\nu |{|^2} + ||\Theta _{k,\alpha }^\nu |{|^2}} \right) \\
  \le & C{\sum\limits_{\nu  = e,i} {} \left( {||\partial _t^k({u^\nu },{\Theta ^\nu })||_{|\alpha | - 1}^2 + ||\partial _t^k{Q^\nu }||_{|\alpha |}^2} \right) +C ||\partial _t^kF||_{|\alpha | - 1}^2}  + C\sum\limits_{\nu  = e,i} {} |||{V^\nu }|||_s^2|||Z||{|_s}. \\
\end{split}\end{equation}
\end{lemma}

\noindent {\bf Proof.} For all $k \in \N$ and $\alpha \in \N^3$ with
$1 \leq k+|\alpha| \leq
   s$, multiplying \eqref{3.3} by  $A_0^\nu(p^\nu, \theta^\nu) V_{k,\alpha }^\nu $
   and taking integrations  in $x$ over $\T$, we obtain
\begin{equation}
\label{3.7}
\begin{split}
\frac{d}{{dt}}\left\langle {A_0^\nu ({p^\nu },{\theta ^\nu
})V_{k,\alpha }^\nu ,V_{k,\alpha }^\nu } \right\rangle  = I_1^\nu  +
I_2^\nu  + I_3^\nu  + I_4^\nu ,\quad \nu = e,i,
  \end{split}
\end{equation}
with
\begin{equation*}
\begin{split}
I_1^\nu  = \left\langle {{\partial _t}A_0^\nu ({p^\nu },{\theta ^\nu
})V_{k,\alpha }^\nu ,V_{k,\alpha }^\nu } \right\rangle , \quad
I_2^\nu  = \left\langle {{B^\nu }({p^\nu },{u^\nu },{\theta ^\nu
},x)V_{k,\alpha }^\nu ,V_{k,\alpha }^\nu } \right\rangle , \quad \nu
= e,i,
  \end{split}
\end{equation*}
%
%
\begin{equation*}
\begin{split}
I_3^\nu  = 2\left\langle {A_0^\nu({p^\nu },{\theta ^\nu }) \partial
_t^k{\partial ^\alpha }{K^\nu },V_{k,\alpha }^\nu } \right\rangle,
 \quad
I_4^\nu  = 2\left\langle {A_0^\nu({p^\nu },{\theta ^\nu })
g_{k,\alpha }^\nu ,V_{k,\alpha }^\nu } \right\rangle,\quad \nu =
e,i.
  \end{split}
\end{equation*}
In the following, we estimate every term on the right-hand side of
\eqref{3.7}.

For the first term $I_1^\nu$, by  \eqref{2.1*}, the third equation
in \eqref{1.1*} and the Sobolev embedding theorem \cite{Ev98}, we
obtain
\begin{equation*}
\begin{split}
||{\partial _t}({p^\nu },{\theta ^\nu })|{|_{{L^\infty }}} \le
C||({u^\nu },\nabla {u^\nu },{\Theta ^\nu })|{|_{{L^\infty }}} \le
C|||{V^\nu }||{|_s}, \quad \nu = e,i.
  \end{split}
\end{equation*}
and then
\begin{equation}
\label{3.8}
\begin{split}
|I_1^\nu | = |\left\langle {{\partial _t}A_0^\nu ({p^\nu },{\theta
^\nu })V_{k,\alpha }^\nu ,V_{k,\alpha }^\nu } \right\rangle | \le
C||{\partial _t}({p^\nu },{\theta ^\nu })|{|_{{L^\infty
}}}||V_{k,\alpha }^\nu |{|^2} \le C|||{V^\nu }|||_s^3 ,\quad \nu =
e,i.
  \end{split}
\end{equation}

For the second term $I_2^\nu$, with the help of \eqref{2.8c},
\eqref{3.1}-\eqref{3.1b} and the fact that the matrix ${{B^\nu
}({p^\nu },{u^\nu },{\theta ^\nu },x)} $ is antisymmetric at the
point $({p^\nu },{u^\nu },{\theta ^\nu }) = ({\bar p^\nu },0,1) $,
we have
\begin{equation}
\label{3.9}
\begin{split}
|I_2^\nu | = |\left\langle {{B^\nu }({p^\nu },{u^\nu },{\theta ^\nu
},x)V_{k,\alpha }^\nu ,V_{k,\alpha }^\nu } \right\rangle | \le
C|||{V^\nu }|||_s^3,\quad \nu = e,i.
  \end{split}
\end{equation}

For the third term $I_3^\nu$, we split it as
\begin{equation}
\label{3.10}
\begin{split}
I_3^\nu  =& 2\left\langle {A_0^\nu ({p^\nu },{\theta ^\nu })\partial
_t^k{\partial ^\alpha }{K^\nu },V_{k,\alpha }^\nu } \right\rangle\\
=& 2\left\langle {A_0^\nu({p^\nu },{\theta ^\nu }) \partial
_t^k{\partial ^\alpha }K_1^\nu ,V_{k,\alpha }^\nu } \right\rangle  +
2\left\langle {A_0^\nu({p^\nu },{\theta ^\nu }) \partial
_t^k{\partial ^\alpha }K_2^\nu ,V_{k,\alpha }^\nu } \right\rangle ,
\quad \nu = e,i.
  \end{split}
\end{equation}
with
\begin{equation}
\label{3.11}
\begin{split}
K_1^\nu  = K_1^\nu \left( {{u^\nu },{\theta ^\nu },F,G,x} \right) =
\left( {\begin{array}{*{20}{c}}
 \vspace{3mm}  { - \displaystyle\frac{1}{{{\theta ^\nu }}}{\Theta ^\nu }}  \\
 \vspace{3mm}     {{q_\nu }\left( {F + {u^\nu } \times \left( {\bar B + G} \right)} \right) - {u^\nu }}  \\
   { - {\Theta ^\nu }}  \\
\end{array}} \right), \quad \nu =
e,i,
  \end{split}
\end{equation}
\begin{equation}
\label{3.11*}
\begin{split}
K_2^\nu  = K_2^\nu \left( {{u^\nu },{\theta ^\nu }} \right) = \left(
{\begin{array}{*{20}{c}}
 \vspace{3mm}   {\displaystyle\frac{1}{{2{\theta ^\nu }}}{{\left| {{u^\nu }} \right|}^2}}  \\
\vspace{3mm}    0  \\
   {\displaystyle\frac{1}{2}{{\left| {{u^\nu }} \right|}^2}}  \\
\end{array}} \right),
 \quad \nu =
e,i.
  \end{split}
\end{equation}
Due to the fact that $K_2^\nu$ is a quadratic term of $u^\nu$, we
get easily
\begin{equation}
\label{3.12}
\begin{split}
\left| {\left\langle {A_0^\nu({p^\nu },{\theta ^\nu }) \partial
_t^k{\partial ^\alpha }K_2^\nu ,V_{k,\alpha }^\nu } \right\rangle }
\right| \le C{\left\| {A_0^\nu }({p^\nu },{\theta ^\nu })
\right\|_{{L^\infty }}}\left\| {\partial _t^k{\partial ^\alpha
}K_2^\nu } \right\|\left\| {V_{k,\alpha }^\nu } \right\| \le
C|||{V^\nu }|||_s^3 ,
 \quad \nu =
e,i.
  \end{split}
\end{equation}
For the first term on the right-hand side of \eqref{3.10}, by
\eqref{2.8} and \eqref{3.11}, we have
\begin{equation}
\label{3.13}
\begin{split}
& 2\left\langle {A_0^\nu({p^\nu },{\theta ^\nu }) \partial
_t^k{\partial ^\alpha }K_1^\nu ,V_{k,\alpha }^\nu } \right\rangle
  - 2{q_\nu }\left\langle {{F_{k,\alpha }},{n^\nu }u_{k,\alpha }^\nu } \right\rangle  \\
  = & - 2\left\langle {{n^\nu }u_{k,\alpha }^\nu ,u_{k,\alpha }^\nu } \right\rangle
  + 2{q_\nu }\left\langle {{{\left( {{u^\nu } \times G} \right)}_{k,\alpha }},{n^\nu }
  u_{k,\alpha }^\nu } \right\rangle  + 2\left\langle {{p^\nu }\left( {
  \frac{{\Theta _{k,\alpha }^\nu }}{{{\theta ^\nu }}} - {{\left( {
  \frac{{{\Theta ^\nu }}}{{{\theta ^\nu }}}} \right)}_{k,\alpha }}}
  \right),Q_{k,\alpha }^\nu } \right\rangle  \\
 & + 2\left\langle {\frac{{{p^\nu }}}{{{\theta ^\nu }}}
 \left( {{{\left( {\frac{{{\Theta ^\nu }}}{{{\theta ^\nu }}}}
 \right)}_{k,\alpha }} - 2\frac{{\Theta _{k,\alpha }^\nu }}{{{\theta ^\nu }}}}
  \right),\Theta _{k,\alpha }^\nu } \right\rangle ,
  \quad \nu =
e,i.
  \end{split}
\end{equation}
For the second term on the right-hand side of \eqref{3.13}, by
\eqref{3.1} and the Leibniz formula, we obtain
\begin{equation}
\label{3.14}
\begin{split}
&
 \left| {\left\langle {{{\left( {{u^\nu } \times G} \right)}_{k,\alpha }},{n^\nu }u_{k,\alpha }^\nu } \right\rangle } \right| \\
  \le & \left| {\left\langle {{u^\nu } \times {G_{k,\alpha }},{n^\nu }u_{k,\alpha }^\nu } \right\rangle } \right| + \sum\limits_{l < k} {} C_k^l\left| {\left\langle {\partial _t^{k - l}{u^\nu } \times {G_{l,\alpha }},{n^\nu }u_{k,\alpha }^\nu } \right\rangle } \right| \\
  & + \sum\limits_{\beta  < \alpha } {} C_\alpha ^\beta \left| {\left\langle {{\partial ^{\alpha  - \beta }}{u^\nu } \times {G_{k,\beta }},{n^\nu }u_{k,\alpha }^\nu } \right\rangle } \right| + \sum\limits_{\beta  < \alpha ,l < k} {} C_\alpha ^\beta C_k^l\left| {\left\langle {u_{k - l,\alpha  - \beta }^\nu  \times {G_{l,\beta }},{n^\nu }u_{k,\alpha }^\nu } \right\rangle } \right| \\
  \le & C|||{V^\nu }|||_s^2|||Z||{|_s}.
 \quad \nu = e,i.
  \end{split}
\end{equation}
We only give the estimate for the most complicate term
$\displaystyle \sum\limits_{\beta  < \alpha ,l < k} {} C_\alpha
^\beta C_k^l\left| {\left\langle {u_{k - l,\alpha  - \beta }^\nu
\times {G_{l,\beta }},{n^\nu }u_{k,\alpha }^\nu } \right\rangle }
\right| $ on the right-hand side of \eqref{3.14}, the estimates for
other terms are omitted here for the sake of simplicity
\begin{equation}
\label{3.14a}
\begin{split}
 \left| {\left\langle {u_{k - l,\alpha  - \beta }^\nu  \times {G_{l,\beta }},{n^\nu }u_{k,\alpha }^\nu } \right\rangle } \right|
  \le &||u_{k - l,\alpha  - \beta }^\nu ||||{G_{l,\beta }}|{|_{{L^\infty }}}||{n^\nu }u_{k,\alpha }^\nu || \\
  \le & C||u_{k - l,\alpha  - \beta }^\nu ||||{G_{l,\beta }}|{|_2}||u_{k,\alpha }^\nu ||\\
   \le & C|||{V^\nu }|||_s^2|||Z|||_s,
 \quad \nu = e,i, \quad \mbox{when }~ l+|\beta| \le 1,
  \end{split}
\end{equation}
and
\begin{equation}
\label{3.14b}
\begin{split}
 \left| {\left\langle {u_{k - l,\alpha  - \beta }^\nu  \times {G_{l,\beta }},{n^\nu }u_{k,\alpha }^\nu }
  \right\rangle } \right| \le & ||u_{k - l,\alpha  - \beta }^\nu |{|_{{L^\infty }}}||{G_{l,\beta }}||||{n^\nu }u_{k,\alpha }^\nu || \\
  \le & C||u_{k - l,\alpha  - \beta }^\nu |{|_2}||{G_{l,\beta }}|||u_{k,\alpha }^\nu || \\
  \le & C|||{V^\nu }|||_s^2|||Z||{|_s},
 \quad \nu = e,i, \quad \mbox{when }~ l+|\beta| \ge 2.
  \end{split}
\end{equation}
Similarly, for the third term on the right-hand side of
\eqref{3.13}, we have
\begin{equation}
\label{3.15}
\begin{split}
\left| {\left\langle {{p^\nu }\left( {\frac{{\Theta _{k,\alpha }^\nu
}}{{{\theta ^\nu }}} - {{\left( {\frac{{{\Theta ^\nu }}}{{{\theta
^\nu }}}} \right)}_{k,\alpha }}} \right),Q_{k,\alpha }^\nu }
\right\rangle } \right| \le C|||{\Theta ^\nu }|||_s^2|||{Q^\nu
}||{|_s},
 \quad \nu = e,i.
  \end{split}
\end{equation}
For the last term on the right-hand side of \eqref{3.13}, it follows
from the Leibniz formula that
\begin{equation}
\label{3.16}
\begin{split}
& \left\langle {\frac{{{p^\nu }}}{{{\theta ^\nu }}}\left( {{{\left(
{\frac{{{\Theta ^\nu }}}{{{\theta ^\nu }}}}
  \right)}_{k,\alpha }} - 2\frac{{\Theta _{k,\alpha }^\nu }}{{{\theta ^\nu }}}} \right),\Theta _{k,\alpha }^\nu } \right\rangle  \\
  = & \left\langle {{n^\nu }\left( {{{\left( {\frac{{{\Theta ^\nu }}}{{{\theta ^\nu }}}} \right)}_{k,\alpha }}
  - 2\frac{{\Theta _{k,\alpha }^\nu }}{{{\theta ^\nu }}}} \right),\Theta _{k,\alpha }^\nu } \right\rangle  \\
  = & - \left\langle {\frac{{{n^\nu }}}{{{\theta ^\nu }}}\Theta _{k,\alpha }^\nu ,\Theta _{k,\alpha }^\nu }
  \right\rangle  + \sum\limits_{l < k} {} C_k^l\left\langle {{n^\nu }\partial _t^{k - l}\left( {\frac{1}{{{\theta ^\nu }}}}
   \right)\Theta _{l,\alpha }^\nu ,\Theta _{k,\alpha }^\nu } \right\rangle  \\
 & + \sum\limits_{\beta  < \alpha } {} C_\alpha ^\beta \left\langle {{n^\nu }{\partial ^{\alpha
   - \beta }}\left( {\frac{1}{{{\theta ^\nu }}}} \right)\Theta _{k,\beta }^\nu ,\Theta _{k,\alpha }^\nu }
    \right\rangle  + \sum\limits_{\beta  < \alpha ,l < k} {} C_\alpha ^\beta C_k^l\left\langle {{n^\nu }
    {{\left( {\frac{1}{{{\theta ^\nu }}}} \right)}_{k - l,\alpha  - \beta }}\Theta _{l,\beta }^\nu ,
    \Theta _{k,\alpha }^\nu } \right\rangle  \\
  \le&  - \left\langle {\frac{{{n^\nu }}}{{{\theta ^\nu }}}\Theta _{k,\alpha }^\nu ,
  \Theta _{k,\alpha }^\nu } \right\rangle  + C|||{\Theta ^\nu }|||_s^3.
 \quad \nu = e,i.
  \end{split}
\end{equation}
Thus from \eqref{3.10}-\eqref{3.16}, for $\nu = e,i$, we obtain
\begin{equation}
\label{3.17}
\begin{split}
 I_3^\nu = & 2\left\langle {A_0^\nu({p^\nu },{\theta ^\nu }) \partial _t^k{\partial ^\alpha
}{K^\nu },V_{k,\alpha }^\nu } \right\rangle\\
\le &
  2{q_\nu }\left\langle
{{F_{k,\alpha }},{n^\nu }u_{k,\alpha }^\nu } \right\rangle  -
2\left\langle {{n^\nu }u_{k,\alpha }^\nu ,u_{k,\alpha }^\nu }
\right\rangle  - 2\left\langle {\frac{{{n^\nu }}}{{{\theta ^\nu
}}}\Theta _{k,\alpha }^\nu ,\Theta _{k,\alpha }^\nu } \right\rangle
+  C|||{V^\nu }|||_s^2|||Z||{|_s}.
  \end{split}
\end{equation}

Next, we begin to estimate $I_4^\nu = 2\left\langle {A_0^\nu({p^\nu
},{\theta ^\nu }) g_{k,\alpha }^\nu ,V_{k,\alpha }^\nu }
\right\rangle$. By \eqref{3.5}, we split $g^\nu_{k,\alpha}$ as
\[g_{k,\alpha }^\nu  = g_{k,\alpha }^{\nu 1} + g_{k,\alpha }^{\nu 2}, \quad \nu = e,i,\]
with
\[g_{k,\alpha }^{\nu 1} = \sum\limits_{j = 1}^3 {} \left( {A_j^\nu
\left( {{u^\nu },{\theta ^\nu }} \right){\partial _j}V_{k,\alpha }^\nu
 - \partial _t^k{\partial ^\alpha }\left( {A_j^\nu \left( {{u^\nu },
 {\theta ^\nu }} \right){\partial _j}{V^\nu }} \right)} \right),
  \quad \nu = e,i,\]
\[g_{k,\alpha }^{\nu 2} = {L^\nu }\left( x \right)V_{k,\alpha }^\nu
 - \partial _t^k{\partial ^\alpha }\left( {{L^\nu }\left( x \right)}
  V^\nu\right), \quad \nu = e,i.\]
We first establish the estimate for $g_{k,\alpha }^{\nu 1}$. Using
\eqref{2.6a} and the Leibniz formula, we have
 \begin{equation*}
\begin{split}
&
 A_j^\nu \left( {{u^\nu },{\theta ^\nu }} \right){\partial _j}V_{k,\alpha }^\nu  - \partial _t^k{\partial ^\alpha }\left( {A_j^\nu \left( {{u^\nu },{\theta ^\nu }} \right){\partial _j}{V^\nu }} \right) \\
  = & \sum\limits_{l < k} {} C_k^l\partial _t^{k - l}u_j^\nu {\partial _j}Q_{l,\alpha }^\nu  + \sum\limits_{\beta  < \alpha } {} C_\alpha ^\beta {\partial ^{\alpha  - \beta }}u_j^\nu {\partial _j}Q_{k,\beta }^\nu  + \sum\limits_{\beta  < \alpha ,l < k} {} C_\alpha ^\beta C_k^l\partial _t^{k - l}{\partial ^{\alpha  - \beta }}u_j^\nu {\partial _j}Q_{l,\beta }^\nu  \\
  +& \sum\limits_{l < k} {} C_k^l\partial _t^{k - l}{\Theta ^\nu }{\partial _j}Q_{l,\alpha }^\nu  + \sum\limits_{\beta  < \alpha } {} C_\alpha ^\beta {\partial ^{\alpha  - \beta }}{\Theta ^\nu }{\partial _j}Q_{k,\beta }^\nu  + \sum\limits_{\beta  < \alpha ,l < k} {} C_\alpha ^\beta C_k^l\Theta _{k - l,\alpha  - \beta }^\nu {\partial _j}Q_{l,\beta }^\nu  \\
  +& \sum\limits_{l < k} {} C_k^l\partial _t^{k - l}u_j^\nu {\partial _j}u_{l,\alpha }^\nu  + \sum\limits_{\beta  < \alpha } {} C_\alpha ^\beta {\partial ^{\alpha  - \beta }}u_j^\nu {\partial _j}u_{k,\beta }^\nu  + \sum\limits_{\beta  < \alpha ,l < k} {} C_\alpha ^\beta C_k^l\partial _t^{k - l}{\partial ^{\alpha  - \beta }}u_j^\nu {\partial _j}u_{l,\beta }^\nu  \\
  +& \sum\limits_{l < k} {} C_k^l\partial _t^{k - l}{\Theta ^\nu }\partial _t^l{\partial ^\alpha }{\partial _j}{u^\nu } + \sum\limits_{\beta  < \alpha } {} C_\alpha ^\beta {\partial ^{\alpha  - \beta }}{\Theta ^\nu }{\partial _j}u_{k,\beta }^\nu  + \sum\limits_{\beta  < \alpha ,l < k} {} C_\alpha ^\beta C_k^l\Theta _{k - l,\alpha  - \beta }^\nu {\partial _j}u_{l,\beta }^\nu  \\
  +& \sum\limits_{l < k} {} C_k^l\partial _t^{k - l}u_j^\nu {\partial _j}\Theta _{l,\alpha }^\nu  + \sum\limits_{\beta  < \alpha } {} C_\alpha ^\beta {\partial ^{\alpha  - \beta }}u_j^\nu {\partial _j}\Theta _{k,\beta }^\nu  + \sum\limits_{\beta  < \alpha ,l < k} {} C_\alpha ^\beta C_k^l\partial _t^{k - l}{\partial ^{\alpha  - \beta }}u_j^\nu {\partial _j}\Theta _{l,\beta }^\nu  \\
  \le& C|||{V^\nu }|||_s^2,
\end{split}
\end{equation*}
which implies
\begin{equation}
\label{3.17a}
\begin{split}
 \left\|g_{k,\alpha }^{\nu 1}\right\| = \left\|\sum\limits_{j = 1}^3 {} \left(
{A_j^\nu \left( {{u^\nu },{\theta ^\nu }} \right){\partial
_j}V_{k,\alpha }^\nu  - \partial _t^k{\partial ^\alpha }\left(
{A_j^\nu \left( {{u^\nu },{\theta ^\nu }} \right){\partial _j}{V^\nu
}} \right)} \right)\right\| \le C|||{V^\nu }|||_s^2,
  \end{split}
\end{equation}
and then
\begin{equation}
\label{3.17b}
\begin{split}
2\left\langle {A_0^\nu({p^\nu },{\theta ^\nu }) g_{k,\alpha }^{\nu
1},V_{k,\alpha }^\nu } \right\rangle  \le C{\left\| {A_0^\nu
}({p^\nu },{\theta ^\nu }) \right\|_{{L^\infty }}}\left\|
{g_{k,\alpha }^{\nu 1}} \right\|\left\| {V_{k,\alpha }^\nu }
\right\| \le C|||{V^\nu }|||_s^3.
  \end{split}
\end{equation}
On the other hand, it follows from \eqref{2.1a}, \eqref{2.6b}  and
\eqref{3.1}-\eqref{3.1b} that matrix function $L(x)$ is regular
bounded in $\T$. 
Then by using the Leibniz formula, \eqref{2.8} and \eqref{2.8b}, we
get
\begin{equation}
\label{3.18*}
\begin{split}
& \left\langle {A_0^\nu({p^\nu },{\theta ^\nu }) g_{k,\alpha }^{\nu 2},V_{k,\alpha }^\nu } \right\rangle\\
   =
 &  \sum\limits_{\beta  < \alpha } {} C_\alpha ^\beta \left\langle {\frac{{{p^\nu }}}{{{\theta ^\nu }}}{\partial ^{\alpha  - \beta }}\left( {\nabla {{\bar q}^\nu }} \right)u_{k,\beta }^\nu ,\Theta _{k,\alpha }^\nu } \right\rangle
  - \sum\limits_{\beta  < \alpha } {} C_\alpha ^\beta \left\langle {{p^\nu }{\partial ^{\alpha  - \beta }}\left( {\nabla {{\bar q}^\nu }} \right)u_{k,\beta }^\nu ,Q_{k,\alpha }^\nu } \right\rangle  \\
   & - \sum\limits_{\beta  < \alpha } {} C_\alpha ^\beta \left\langle {\frac{{{p^\nu }}}{{{\theta ^\nu }}}{\partial ^{\alpha  - \beta }}\left( {\nabla {{\bar q}^\nu }} \right)\Theta _{k,\beta }^\nu ,u_{k,\alpha }^\nu } \right\rangle , \\
   \end{split}
\end{equation}
which implies
\begin{equation}
\label{3.18}
\begin{split}
2\left| {\left\langle {A_0^\nu ({p^\nu },{\theta ^\nu })g_{k,\alpha
}^{\nu 2},V_{k,\alpha }^\nu } \right\rangle } \right| \le C\left\|
{(\partial _t^k{u^\nu },\partial _t^k{\Theta ^\nu })}
\right\|_{\left| \alpha  \right| - 1}^2 + \varepsilon {\left\|
{(u_{k,\alpha }^\nu ,\Theta _{k,\alpha }^\nu )} \right\|^2} +
C\left\| {\partial _t^k{Q^\nu }} \right\|_{\left| \alpha
\right|}^2, \end{split}
\end{equation}
and thus, 
\begin{equation}
\label{3.19}
\begin{split}
 \left| {I_4^\nu } \right| \le & 2\left| {\left\langle {A_0^\nu ({p^\nu }
 ,{\theta ^\nu })g_{k,\alpha }^{\nu 1},V_{k,\alpha }^\nu } \right\rangle } \right|
  + 2\left| {\left\langle {A_0^\nu ({p^\nu },{\theta ^\nu })g_{k,\alpha }^{\nu 2},
  V_{k,\alpha }^\nu } \right\rangle } \right| \\
  \le & C\left\| {(\partial _t^k{u^\nu },\partial _t^k{\Theta ^\nu })}
   \right\|_{\left| \alpha  \right| - 1}^2 + \varepsilon {\left\|
   {(u_{k,\alpha }^\nu ,\Theta _{k,\alpha }^\nu )} \right\|^2} + C\left\|
   {\partial _t^k{Q^\nu }} \right\|_{\left| \alpha  \right|}^2 + C|||{V^\nu }|||_s^3,\quad \nu = e,i.
\end{split}
\end{equation}
Moreover, a normal energy estimate for (\ref{3.4})
gives %
\begin{equation}
\label{3.20}
\begin{split}
\frac{d}{{dt}}\left( {{{\left\| {{F_{k,\alpha }}} \right\|}^2} +
{{\left\| {{G_{k,\alpha }}} \right\|}^2}} \right) + 2\left\langle
{{{\left( {{n^i}{u^i} - {n^e}{u^e}} \right)}_{k,\alpha
}},{F_{k,\alpha }}} \right\rangle  = 0.
\end{split}
\end{equation}
Hence, combining \eqref{3.8}-\eqref{3.9}, \eqref{3.17} and
\eqref{3.19}-\eqref{3.20}, we obtain
\begin{equation}
\label{3.21}
\begin{split}
&  \frac{d}{{dt}}\left( {\sum\limits_{\nu  = e,i} {} \left\langle
{A_0^\nu ({p^\nu },{\theta ^\nu })V_{k,\alpha }^\nu ,V_{k,\alpha
}^\nu } \right\rangle  + {{\left\| {{F_{k,\alpha }}} \right\|}^2} +
{{\left\| {{G_{k,\alpha }}} \right\|}^2}} \right) +
2\sum\limits_{\nu  = e,i} {} \left\langle {{n^\nu }u_{k,\alpha }^\nu
,u_{k,\alpha }^\nu } \right\rangle
\\
 & + 2\sum\limits_{\nu  = e,i} {} \left\langle {\frac{{{n^\nu }}}{{{\theta ^\nu }}}\Theta _{k,\alpha }^\nu ,\Theta _{k,\alpha }^\nu } \right\rangle  \\
  \le & \varepsilon \sum\limits_{\nu  = e,i} {} {\left\| {(u_{k,\alpha }^\nu ,\Theta _{k,\alpha }^\nu )} \right\|^2} + C\sum\limits_{\nu  = e,i} {} \left\| {(\partial _t^k{u^\nu },\partial _t^k{\Theta ^\nu })} \right\|_{\left| \alpha  \right| - 1}^2 + C\sum\limits_{\nu  = e,i} {} \left\| {\partial _t^k{Q^\nu }} \right\|_{\left| \alpha  \right|}^2 \\
  & + 2\left\langle {{{\left( {{n^e}{u^e} - {n^i}{u^i}} \right)}_{k,\alpha }} - \left( {{n^e}u_{k,\alpha }^e - {n^i}u_{k,\alpha }^i} \right),{F_{k,\alpha }}} \right\rangle  + C\sum\limits_{\nu  = e,i} {} |||{V^\nu }|||_s^2|||Z||{|_s}. \\
\end{split}
\end{equation}

Next, let us estimate the term $2\left\langle {{{\left( {{n^e}{u^e}
- {n^i}{u^i}} \right)}_{k,\alpha }} - \left( {{n^e}u_{k,\alpha }^e -
{n^i}u_{k,\alpha }^i} \right),{F_{k,\alpha }}} \right\rangle$ on the
right-hand side of \eqref{3.21}. With the help of $|\alpha| \ge 1$,
\eqref{2.4} and an integration by parts, we have
\begin{equation}
\label{3.22}
\begin{split}
& \left| {\left\langle {{{\left( {{n^e}{u^e} - {n^i}{u^i}} \right)}_{k,\alpha }} - \left( {{n^e}u_{k,\alpha }^e - {n^i}u_{k,\alpha }^i} \right),{F_{k,\alpha }}} \right\rangle } \right| \\
  \le & \left| {\left\langle {{{\left( {{{\bar n}^e}{u^e}} \right)}_{k,\alpha }} - {{\bar n}^e}u_{k,\alpha }^e,{F_{k,\alpha }}} \right\rangle } \right| + \left| {\left\langle {{{\left( {{N^e}{u^e}} \right)}_{k,\alpha }} - {N^e}u_{k,\alpha }^e,{F_{k,\alpha }}} \right\rangle } \right| \\
  & + \left| {\left\langle {{{\left( {{{\bar n}^i}{u^i}} \right)}_{k,\alpha }} - {{\bar n}^i}u_{k,\alpha }^i,{F_{k,\alpha }}} \right\rangle } \right| + \left| {\left\langle {{{\left( {{N^i}{u^i}} \right)}_{k,\alpha }} - {N^i}u_{k,\alpha }^i,{F_{k,\alpha }}} \right\rangle } \right| \\
  \le & \left| {\left\langle {{\partial _x}\left( {{{\left( {{{\bar n}^e}{u^e}} \right)}_{k,\alpha }} - {{\bar n}^e}u_{k,\alpha }^e} \right),{F_{k,{\alpha _*}}}} \right\rangle } \right| + \left| {\left\langle {{\partial _x}\left( {{{\left( {{{\bar n}^i}{u^i}} \right)}_{k,\alpha }} - {{\bar n}^i}u_{k,\alpha }^i} \right),{F_{k,{\alpha _*}}}} \right\rangle } \right| \\
  & + C\sum\limits_{\nu  = e,i} {} |||{V^\nu }|||_s^2|||Z||{|_s} \\
  \le & C{\left\| {\partial _t^kF} \right\|_{\left| \alpha  \right| - 1}}\sum\limits_{\nu  = e,i} {} \left( {\left\| {u_{k,\alpha }^\nu } \right\| + {{\left\| {\partial _t^k{u^\nu }} \right\|}_{\left| \alpha  \right| - 1}}} \right) + C\sum\limits_{\nu  = e,i} {} |||{V^\nu }|||_s^2|||Z||{|_s} \\
  \le & \varepsilon \sum\limits_{\nu  = e,i} {} {\left\| {u_{k,\alpha }^\nu } \right\|^2} + C\left\| {\partial _t^kF} \right\|_{\left| \alpha  \right| - 1}^2 + C\sum\limits_{\nu  = e,i} {} \left\| {\partial _t^k{u^\nu }} \right\|_{\left| \alpha  \right| - 1}^2 + C\sum\limits_{\nu  = e,i} {} |||{V^\nu }|||_s^2|||Z||{|_s}, \\
\end{split}
\end{equation}
where $\alpha_* \in \mathbb{N}^3$ with $|\alpha_*| = |\alpha| - 1$
for $ |\alpha| \ge 1$.

Then, taking $\varepsilon >0$ sufficiently small (for example
$\varepsilon = \displaystyle\frac{1}{17}$ ), by
\eqref{3.1}-\eqref{3.1b}, the combination of
\eqref{3.21}-\eqref{3.22} yields \eqref{3.6}.\hfill $\Box$

\begin{remark} Lemma \ref{L3.3} is valid for $|\alpha| \geq 1$. The
following Lemma concerns the $L^2$ estimates for $\pa_t^k V^\nu$
(i.e. $\alpha=0$), which is a starting point for employing the
argument of induction.
\end{remark}

\begin{lemma}
\label{L3.4} 
Assume that the conditions of Theorem \ref{T1} holds and $\omega_T$
is sufficiently small independent of $T$, then there exists a
positive constant $C_0$ such that, for all $0 \leq k \leq s$,
  it holds
\begin{equation}
\label{3.25}
\begin{split}
&  \frac{d}{{dt}}  \left( {\sum\limits_{\nu  = e,i} {} \left\langle
{A_0^\nu({p^\nu },{\theta ^\nu }) V_{k,0}^\nu ,V_{k,0}^\nu }
\right\rangle  + ||{F_{k,0}}|{|^2} + ||{G_{k,0}}|{|^2}} \right) +
{C_0}\sum\limits_{\nu  = e,i} {} \left(
{||u_{k,0}^\nu |{|^2} + ||\Theta _{k,0}^\nu |{|^2}} \right)\\
 \le  &
C\sum\limits_{\nu  = e,i} {} |||{V^\nu }|||_s^2|||Z||{|_s}.
\end{split}\end{equation}
\end{lemma}
\noindent{\textbf{Proof.}} Carefully checking the procedures of the
proof for Lemma \ref{L3.3}, we shall prove
\begin{equation}
\label{3.26}
\begin{split}
2\left| {\left\langle {A_0^\nu ({p^\nu },{\theta ^\nu })g_{k,0}^{\nu
2},V_{k,0}^\nu } \right\rangle } \right| \le C|||{V^\nu
}|||_s^2|||Z||{|_s},\end{split}\end{equation}
and
\begin{equation}
\label{3.27}
\begin{split}
\left| {\left\langle {{{\left( {{n^e}{u^e} - {n^i}{u^i}}
\right)}_{k,\alpha }} - \left( {{n^e}u_{k,\alpha }^e -
{n^i}u_{k,\alpha }^i} \right),{F_{k,\alpha }}} \right\rangle }
\right| \le C\sum\limits_{\nu  = e,i} {} |||{V^\nu
}|||_s^2|||Z||{|_s},\end{split}\end{equation}
which are correspondence to \eqref{3.18} and \eqref{3.22},
respectively. Obviously, \eqref{3.26} can be easily obtained through
the following calculation
\[g_{k,0}^{\nu 2} = {L^\nu }(x)V_{k,0}^\nu  - \partial _t^k\left( {{L^\nu }(x){V^\nu }} \right) = 0, \quad \nu =e, i.\]
Next, \eqref{3.27} follows from the following calculations
\begin{equation*}
\begin{split}
& \left| {\left\langle {{{\left( {{n^e}{u^e} - {n^i}{u^i}}
\right)}_{k,\alpha }} - \left( {{n^e}u_{k,\alpha }^e -
{n^i}u_{k,\alpha }^i} \right),{F_{k,\alpha }}} \right\rangle }
\right|\\
 = & \left| {\left\langle {{{\left( {{N^e}{u^e} - {N^i}{u^i}}
\right)}_{k,\alpha }} - \left( {{N^e}u_{k,\alpha }^e -
{N^i}u_{k,\alpha }^i} \right),{F_{k,\alpha }}} \right\rangle }
\right| \le C\sum\limits_{\nu  = e,i} {} |||{V^\nu
}|||_s^2|||Z||{|_s}.\end{split}\end{equation*}
 \hfill $\Box$

\vspace{3mm}

\subsection{Recurrence relationship}%
%
In order to prove Theorem \ref{T1}, we have to control the terms
$\left\| {\partial _t^kF} \right\|_{\left| \alpha  \right| - 1}^2$
and $\left\| {\partial _t^k{Q^\nu }} \right\|_{\left| \alpha
\right|}^2$ appearing on the right-hand side of \eqref{3.6}. This
will be achieved in the following Lemma.
\begin{lemma}
\label{L3.5}
Assume that the conditions of Theorem \ref{T1} holds and $\omega_T$
is sufficiently small independent of $T$, then for all $k \in \N$
and $\alpha \in \N^3$ with $|\alpha| \geq 1$ and $k+|\alpha| \leq
s$, it holds
\begin{equation}
\label{3.28}
\begin{split}
\left\| {\partial _t^k{N^\nu }} \right\|_{|\alpha |}^2 \le C\left(
{\left\| {\partial _t^k{Q^\nu }} \right\|_{|\alpha | - 1}^2 +
\left\| {\partial _t^k{\Theta ^\nu }} \right\|_{|\alpha | - 1}^2}
\right) + C|||{V^\nu }|||_s^2|||Z||{|_s},\quad \nu = e, i,
\end{split}
\end{equation}
\begin{equation}
\label{3.29}
\small\begin{split} \sum\limits_{\nu  = e,i} {} \left\| {\partial
_t^k{Q^\nu }} \right\|_{|\alpha |}^2 \le C\sum\limits_{\nu  = e,i}
{} \left( {\left\| {\partial _t^k\left( {{Q^\nu },{u^\nu },{\Theta
^\nu }} \right)} \right\|_{|\alpha | - 1}^2 + \left\| {\partial
_t^{k + 1}{u^\nu }} \right\|_{|\alpha | - 1}^2} \right) +
C\sum\limits_{\nu = e,i} {} |||{V^\nu }|||_s^2|||Z||{|_s},
\end{split}
\end{equation}
and
\begin{equation}
\label{3.30}
\small\begin{split} \sum\limits_{\nu  = e,i} {} \left\| {\partial
_t^k{F }} \right\|_{|\alpha | - 1}^2 \le C\sum\limits_{\nu  = e,i}
{} \left( {\left\| {\partial _t^k\left( {{Q^\nu },{u^\nu },{\Theta
^\nu }} \right)} \right\|_{|\alpha | - 1}^2 + \left\| {\partial
_t^{k + 1}{u^\nu }} \right\|_{|\alpha | - 1}^2} \right) +
C\sum\limits_{\nu = e,i} {} |||{V^\nu }|||_s^2|||Z||{|_s},
\end{split}
\end{equation}
\end{lemma}

\noindent{\textbf{Proof.}} The first estimate \eqref{3.28} follows
from \eqref{2.4}. 
%
%
Next we prove \eqref{3.29}, rewriting the second equation in
\eqref{2.3} as
\begin{equation}
\label{3.31}
\begin{split}
\nabla {Q^\nu } = {q_\nu }F - \left( {{\partial _t}{u^\nu } + {u^\nu }} \right)
 + {q_\nu }{u^\nu } \times \bar B - \nabla {{\bar q}^\nu }{\Theta ^\nu } + {r^\nu },
\end{split}
\end{equation}
with
\[{r^\nu } =
  - \left( {{u^\nu } \cdot \nabla } \right){u^\nu }
  - {\Theta ^\nu }\nabla {Q^\nu } + {q_\nu }{u^\nu } \times G, \quad \nu = e, i.\]

For $k \in \N$ and $\beta \in \N^3$ with $\beta < \alpha$  and
$k+|\alpha| \leq s$, applying $\pa_t^k \pa^\beta$ to (\ref{3.31})
and taking the inner product with $\left(\nabla Q^\nu\right)  _{k,
\beta}$ in $L^2(\mathbb{T})$, we have
\begin{equation}
\label{3.32}
\begin{split}
 {\left\| {{{\left( {\nabla {Q^\nu }} \right)}_{k,\beta }}} \right\|^2} = & {q_\nu }\left\langle
  {{F_{k,\beta }},{{\left( {\nabla {Q^\nu }} \right)}_{k,\beta }}} \right\rangle
  + \left\langle {u_{k + 1,\beta }^\nu  + u_{k,\beta }^\nu ,{{\left( {\nabla {Q^\nu }}
   \right)}_{k,\beta }}} \right\rangle  \\
  & + {q_\nu }\left\langle {u_{k,\beta }^\nu  \times \bar B,{{\left( {\nabla {Q^\nu }}
   \right)}_{k,\beta }}} \right\rangle  - \left\langle {{{\left( {\nabla {{\bar q}^\nu }
   {\Theta ^\nu }} \right)}_{k,\beta }},{{\left( {\nabla {Q^\nu }} \right)}_{k,\beta }}} \right\rangle
  + \left\langle {r_{k,\beta }^\nu ,{{\left( {\nabla {Q^\nu }} \right)}_{k,\beta }}}
  \right\rangle.
\end{split}
\end{equation}
By the compatibility condition $\nabla \cdot F = N^i -N^e$, we get
\[\begin{array}{l}
 \displaystyle\sum\limits_{\nu  = e,i} {} {q_\nu }\left\langle {{F_{k,\beta }},{{\left( {\nabla {Q^\nu }} \right)}_{k,\beta }}} \right\rangle  = \left\langle {{F_{k,\beta }},{{\left( {\nabla {Q^i} - \nabla {Q^e}} \right)}_{k,\beta }}} \right\rangle  \\
  =  - \left\langle {\nabla  \cdot {F_{k,\beta }},{{\left( {{Q^i} - {Q^e}} \right)}_{k,\beta }}} \right\rangle  = \left\langle {{{\left( {{N^e} - {N^i}} \right)}_{k,\beta }},{{\left( {{Q^i} - {Q^e}} \right)}_{k,\beta }}} \right\rangle  \\
 \end{array}\]

Then \eqref{3.28} together with the Young inequality implies that
\begin{equation}
\label{3.33}
\begin{split}
\left| {\sum\limits_{\nu  = e,i} {} {q_\nu }\left\langle
{{F_{k,\beta }},{{\left( {\nabla {Q^\nu }} \right)}_{k,\beta }}}
\right\rangle } \right| \le C\sum\limits_{\nu  = e,i} {} \left\|
{\partial _t^k\left( {{Q^\nu },{\Theta ^\nu }} \right)}
\right\|_{|\alpha | - 1}^2 + C\sum\limits_{\nu  = e,i} {} |||{V^\nu
}|||_s^2|||Z||{|_s}.
\end{split}
\end{equation}
In a similar way, by the Young inequality, we obtain
\begin{equation}
\label{3.34}
\begin{split}
& \left| {\left\langle {u_{k + 1,\beta }^\nu  + u_{k,\beta }^\nu ,{{\left( {\nabla {Q^\nu }} \right)}_{k,\beta }}} \right\rangle  + {q_\nu }\left\langle {u_{k,\beta }^\nu  \times \bar B,{{\left( {\nabla {Q^\nu }} \right)}_{k,\beta }}} \right\rangle } \right| \\
  \le& C\left( {\left\| {\partial _t^k{u^\nu }} \right\|_{|\alpha | - 1}^2 + \left\| {\partial _t^{k + 1}{u^\nu }} \right\|_{|\alpha | - 1}^2} \right) + \varepsilon {\left\| {{{\left( {\nabla {Q^\nu }} \right)}_{k,\beta }}} \right\|^2}, \\
\end{split}
\end{equation}
and
\begin{equation}
\label{3.35}
\begin{split}
|\left\langle {{{\left( {\nabla {{\bar q}^\nu }{\Theta ^\nu }}
\right)}_{k,\beta }},{{\left( {\nabla {Q^\nu }} \right)}_{k,\beta
}}} \right\rangle|  \le C\left\| {\partial _t^k{\Theta ^\nu }}
\right\|_{|\alpha | - 1}^2 + \varepsilon {\left\| {{{\left( {\nabla
{Q^\nu }} \right)}_{k,\beta }}} \right\|^2}.
\end{split}
\end{equation}
On the other hand, by the Leibniz formula, we easily get%
\begin{equation}
\label{3.36}
\begin{split}
\left| {\left\langle {r_{k,\beta }^\nu ,{{\left( {\nabla {Q^\nu }}
\right)}_{k,\beta }}} \right\rangle } \right| \le C|||{V^\nu
}|||_s^2|||Z||{|_s}, \quad \nu = e, i.
\end{split}
\end{equation}

Taking $\varepsilon >0$ sufficiently small ( for instance,
$\varepsilon = \displaystyle\frac{1}{6}$), and combining
\eqref{3.32}-\eqref{3.36},  we have
\[\sum\limits_{\nu  = e,i} {} {\left\| {{{\left( {\nabla {Q^\nu }} \right)}_{k,\beta }}} \right\|^2} \le C\sum\limits_{\nu  = e,i} {} \left( {\left\| {\partial _t^k\left( {{Q^\nu },{u^\nu },{\Theta ^\nu }} \right)} \right\|_{|\alpha | - 1}^2 + \left\| {\partial _t^{k + 1}{u^\nu }} \right\|_{|\alpha | - 1}^2} \right) + C\sum\limits_{\nu  = e,i} {} |||{V^\nu }|||_s^2|||Z||{|_s}.\]
Then  \eqref{3.29} follows from the summation of these inequalities
for all indexes $\beta < \alpha$.

In the end, by \eqref{3.31}, we have
$$
{q_\nu }F = \nabla {Q^\nu } + \left( {{\partial _t}{u^\nu } + {u^\nu
}} \right) - {q_\nu }{u^\nu } \times \bar B + \nabla {\bar q^\nu
}{\Theta ^\nu } - {r^\nu },$$
this equality  together \eqref{3.29} imply \eqref{3.30}.
 \hfill $\Box$

\vspace{3mm}

Now we give a dissipation estimate for $\left\| {\partial _t^s{Q^\nu
}} \right\|^2$ and a refined estimate of \eqref{3.29} for $ \left\|
{\partial _t^k{Q^\nu }} \right\|^2 $, with $k \le s-1$. It should be
pointed out that the two estimates are necessary, i.e,
 the process by induction in the proof of Theorem \ref{T1} may not
 be closed without of them (see \eqref{3.46}).
\begin{lemma}
\label{L3.6} 
 Assume that the conditions of Theorem \ref{T1} holds and $\omega_T$
is sufficiently small independent of $T$, then for all $k \in
\mathbb{N}$ with $ k \leq s-1$, we have
\begin{equation}
\label{3.37}
\left\| {\partial _t^k{Q^\nu }} \right\|^2 \le C\sum\limits_{\nu  =
e,i} {} \left( {{{\left\| {\partial _t^k\left( {{u^\nu },{\Theta
^\nu }} \right)} \right\|}^2} + {{\left\| {\partial _t^{k + 1}{u^\nu
}} \right\|}^2}} \right) + C\sum\limits_{\nu  = e,i} {} |||{V^\nu
}|||_s^2|||Z||{|_s}, \quad \nu=e,i,
\end{equation}
and
\begin{equation}
\label{3.38}
\begin{split}
{\left\| {\partial _t^s{Q^\nu }} \right\|^2} \le C\left\| {\partial
_t^{s - 1}\left( {{u^\nu },{\Theta ^\nu }} \right)} \right\|_1^2 +
C|||{V^\nu }|||_s^2|||Z||{|_s}, \quad \nu=e,i.
 \end{split}\end{equation}

\end{lemma}
\noindent{\textbf{Proof.}} For $k \in \N$ with $k  \leq s-1$,
applying $\pa_t^k  $ to (\ref{3.31}), we get
\begin{equation}
\label{3.38a}
\begin{split}
\nabla \partial _t^k{Q^\nu } - {q_\nu }\partial _t^kF = {q_\nu
}\partial _t^k{u^\nu } \times \bar B - \partial _t^k{u^\nu } -
\partial _t^{k + 1}{u^\nu } - \nabla {\bar q^\nu }\partial
_t^k{\Theta ^\nu } + \partial _t^k{r^\nu },\quad \nu=e,i.
\end{split}
\end{equation}
Now, we define a potential function $\nabla \psi$ as
$$\nabla \psi  = \bar E - E =  - F, \quad \int_{\T}\psi(t,x)dx =0.$$
Then
$$\nabla  \cdot \left( {{\partial_t^k}F + \nabla {\partial_t^k}\psi} \right) = 0,
\quad \forall~~ 0 \le k \le s-1.$$
From \eqref{3.38a}, we have
\begin{equation}
\label{3.38b}
\begin{split}
\nabla \eta _k^\nu  - {q_\nu }\left( {\partial _t^kF + \nabla
\partial _t^k\psi } \right) = {q_\nu }\partial _t^k{u^\nu } \times
\bar B - \partial _t^k{u^\nu } - \partial _t^{k + 1}{u^\nu } -
\nabla {{\bar q}^\nu }\partial _t^k{\Theta ^\nu } + \partial
_t^k{r^\nu },
 \end{split}\end{equation}
where 
\begin{equation}
\label{3.38c}
\begin{split}
\eta _k^\nu=\partial _t^k \eta ^\nu = \partial _t^k{Q^\nu } + {q_\nu
}\partial _t^k\psi, \quad \nu = e,i.\end{split}
\end{equation}
Due to the fact that $$\left\langle {\partial _t^kF + \nabla
\partial _t^k\psi ,\nabla \eta _k^\nu } \right\rangle  =  -
\left\langle {\partial _t^k\nabla  \cdot F + \Delta \partial
_t^k\psi ,\eta _k^\nu } \right\rangle  = 0,
$$
 we obtain
\begin{equation}
\label{3.38d}
\begin{split}
{\left\| {\nabla \eta _k^\nu } \right\|^2} \le {\left\| {\partial
_t^k{u^\nu }} \right\|^2} + {\left\| {\partial _t^{k + 1}{u^\nu }}
\right\|^2} + C{\left\| {\partial _t^k{\Theta ^\nu }} \right\|^2} +
{\left\| {\partial _t^k{r^\nu }} \right\|^2}.
 \end{split}
\end{equation}
%
%
From \eqref{3.38c}, we have
$$
\partial _t^k{Q^\nu } = \eta _k^\nu  - {q_\nu }\partial _t^k\psi.
$$
For $k=0$, we get
\[{Q^\nu } = {\eta ^\nu } - {q_\nu }\psi , \quad {q^\nu } = {Q^\nu } + {\bar q^\nu }.\]
Since
$${N^\nu } = {n^\nu } - {\bar n^\nu } = \frac{{{p^\nu }}}{{{\theta
^\nu }}} - \frac{{{{\bar p}^\nu }}}{{{{\bar \theta }^\nu }}} =
\frac{{{e^{{q^\nu }}}}}{{{\theta ^\nu }}} - \frac{{{e^{{{\bar q}^\nu
}}}}}{{{{\bar \theta }^\nu }}} = \frac{{{e^{{{\bar q}^\nu
}}}}}{{{\theta ^\nu }}}\left( {{e^{{Q^\nu }}} - 1 - {\Theta ^\nu }}
\right) = \frac{{{e^{{{\bar q}^\nu }}}}}{{{\theta ^\nu }}}\left(
{{e^{{\eta ^\nu }}}{e^{ - {q_\nu }\psi }} - 1 - {\Theta ^\nu }}
\right),$$
we have
$$- \Delta \psi  = {N^i} - {N^e} = \frac{{{e^{{{\bar
q}^i}}}}}{{{\theta ^i}}}\left( {{e^{{\eta ^i}}}{e^{ - \psi }} - 1}
\right) - \frac{{{e^{{{\bar q}^e}}}}}{{{\theta ^e}}}\left(
{{e^{{\eta ^e}}}{e^\psi } - 1} \right) - \frac{{{e^{{{\bar
q}^i}}}}}{{{\theta ^i}}}{\Theta ^i} + \frac{{{e^{{{\bar
q}^e}}}}}{{{\theta ^e}}}{\Theta ^e}.
$$
Thus, $$  - \Delta \psi  + \left( {\frac{{{e^{{{\bar
q}^e}}}}}{{{\theta ^e}}}{e^{{\eta ^e}}} + \frac{{{e^{{{\bar
q}^i}}}}}{{{\theta ^i}}}{e^{{\eta ^i}}}} \right)\psi  =
\frac{{{e^{{{\bar q}^e}}}}}{{{\theta ^e}}}\left( {{e^{{\eta
^e}}}\left( {\psi  - {e^\psi }} \right) + 1} \right) +
\frac{{{e^{{{\bar q}^i}}}}}{{{\theta ^i}}}\left( {{e^{{\eta
^i}}}\left( {\psi  + {e^{ - \psi }}} \right) - 1} \right) +
\frac{{{e^{{{\bar q}^e}}}}}{{{\theta ^e}}}{\Theta ^e} -
\frac{{{e^{{{\bar q}^i}}}}}{{{\theta ^i}}}{\Theta ^i}.
$$
Since $ \left( {\displaystyle\frac{{{e^{{{\bar q}^e}}}}}{{{\theta
^e}}}{e^{{\eta ^e}}} + \frac{{{e^{{{\bar q}^i}}}}}{{{\theta
^i}}}{e^{{\eta ^i}}}} \right)
 \geq \mbox{const.} >0$, taking the inner product
of the previous equality with $\psi $ in $L^2(\T)$ and using an
integration by parts, we get
\begin{equation}
\label{3.38e}
\begin{split}
{\left\| {\nabla \psi } \right\|^2} + {C_0}{\left\| { \psi }
\right\|^2} \le \sum\limits_{\nu  = e,i} {} \left( {{{\left\| {{\eta
^\nu }} \right\|}^2} + {{\left\| {{\Theta ^\nu }} \right\|}^2}}
\right) \le \sum\limits_{\nu  = e,i} {} \left( {{{\left\| {\nabla
{\eta ^\nu }} \right\|}^2} + {{\left\| {{\Theta ^\nu }}
\right\|}^2}} \right),
 \end{split}\end{equation}
where we have used Lemma  \ref{L2.1}.

For $k\ge 1$, due to the fact that
\[\partial _t^k{N^\nu } =  - {q_\nu }\frac{{{e^{{{\bar q}^\nu }}}}}{{{\theta ^\nu }}}{e^{{\eta ^\nu }}}\partial _t^k\psi  + \left( {\partial _t^k\left( {\frac{{{e^{{{\bar q}^\nu }}}}}{{{\theta ^\nu }}}\left( {{e^{{\eta ^\nu }}}{e^{ - {q_\nu }\psi }} - 1} \right)} \right) + {q_\nu }\frac{{{e^{{{\bar q}^\nu }}}}}{{{\theta ^\nu }}}{e^{{\eta ^\nu }}}\partial _t^k\psi } \right) - \partial _t^k\left( {\frac{{{e^{{{\bar q}^\nu }}}}}{{{\theta ^\nu }}}{\Theta ^\nu }} \right),\]
we have
\begin{equation*}
\begin{split}
 & - \Delta \partial _t^k\psi  = \partial _t^k{N^i} - \partial _t^k{N^e} \\
  =  & - \frac{{{e^{{{\bar q}^i}}}}}{{{\theta ^i}}}{e^{{\eta ^i}}}\partial _t^k\psi  + \left( {\partial _t^k\left( {\frac{{{e^{{{\bar q}^i}}}}}{{{\theta ^i}}}\left( {{e^{{\eta ^i}}}{e^{ - \psi }} - 1} \right)} \right) + \frac{{{e^{{{\bar q}^i}}}}}{{{\theta ^i}}}{e^{{\eta ^i}}}\partial _t^k\psi } \right) - \partial _t^k\left( {\frac{{{e^{{{\bar q}^i}}}}}{{{\theta ^i}}}{\Theta ^i}} \right) \\
  & - \frac{{{e^{{{\bar q}^e}}}}}{{{\theta ^e}}}{e^{{\eta ^e}}}\partial _t^k\psi
  - \left( {\partial _t^k\left( {\frac{{{e^{{{\bar q}^e}}}}}{{{\theta ^e}}}\left(
  {{e^{{\eta ^e}}}{e^\psi } - 1} \right)} \right) - \frac{{{e^{{{\bar q}^e}}}}}{{{
  \theta ^e}}}{e^{{\eta ^e}}}\partial _t^k\psi } \right) + \partial _t^k\left(
  {\frac{{{e^{{{\bar q}^e}}}}}{{{\theta ^e}}}{\Theta ^e}} \right), \end{split}
\end{equation*}
which implies that
\begin{equation*}
\begin{split}
 & - \Delta \partial _t^k\psi  + \left( {\frac{{{e^{{{\bar q}^i}}}}}{{{\theta ^i}}}{e^{{\eta ^i}}}
 + \frac{{{e^{{{\bar q}^e}}}}}{{{\theta ^e}}}{e^{{\eta ^e}}}} \right)\partial
 _t^k\psi\\
  = & \left( {\partial _t^k\left( {\frac{{{e^{{{\bar q}^i}}}}}{{{\theta ^i}}}\left( {{e^{{\eta ^i}}}{e^{ - \psi }} - 1} \right)}
   \right) + \frac{{{e^{{{\bar q}^i}}}}}{{{\theta ^i}}}{e^{{\eta ^i}}}\partial _t^k\psi } \right)
   - \partial _t^k\left( {\frac{{{e^{{{\bar q}^i}}}}}{{{\theta ^i}}}{\Theta ^i}} \right) \\
  &  - \left( {\partial _t^k\left( {\frac{{{e^{{{\bar q}^e}}}}}{{{\theta ^e}}}
  \left( {{e^{{\eta ^e}}}{e^\psi } - 1} \right)} \right) - \frac{{{e^{{{\bar q}^e}}}}}
  {{{\theta ^e}}}{e^{{\eta ^e}}}\partial _t^k\psi } \right) + \partial _t^k
  \left( {\frac{{{e^{{{\bar q}^e}}}}}{{{\theta ^e}}}{\Theta ^e}} \right).
\end{split}
\end{equation*}
Since $ \left( {\displaystyle\frac{{{e^{{{\bar q}^e}}}}}{{{\theta
^e}}}{e^{{\eta ^e}}} + \frac{{{e^{{{\bar q}^i}}}}}{{{\theta
^i}}}{e^{{\eta ^i}}}} \right)
 \geq \mbox{const.} >0$, taking the inner product
of the previous equality with $\partial _t^k\psi $ in $L^2(\T)$, by
the Leibniz formula and Lemma  \ref{L2.1}, using an integration by
parts, we get
\begin{equation}
\label{3.38f}
\begin{split}
{\left\| {\nabla \partial _t^k\psi } \right\|^2} + {C_0}{\left\|
{\partial _t^k\psi } \right\|^2} \le \sum\limits_{\nu  = e,i} {}
\left( {{{\left\| \eta _k^\nu \right\|}^2} + {{\left\| {\partial
_t^k{\Theta ^\nu }} \right\|}^2}} \right) \le \sum\limits_{\nu  =
e,i} {} \left( {{{\left\| {\nabla \eta _k^\nu} \right\|}^2} +
{{\left\| {\partial _t^k{\Theta ^\nu }} \right\|}^2}} \right).
 \end{split}\end{equation}
From \eqref{3.38c}-\eqref{3.38f}, we have
\begin{equation*}\begin{split}
& {\left\| {\partial _t^k{Q^\nu }} \right\|^2} \le {\left\| {\eta _k^\nu } \right\|^2} + {\left\| {\partial _t^k\psi } \right\|^2} \\
  \le & C\sum\limits_{\nu  = e,i} {} \left( {{{\left\| {\partial _t^k{u^\nu }} \right\|}^2} + {{\left\| {\partial _t^{k + 1}{u^\nu }} \right\|}^2} + C{{\left\| {\partial _t^k{\Theta ^\nu }} \right\|}^2} + {{\left\| {\partial _t^k{r^\nu }} \right\|}^2}} \right) \\
  \le & C\sum\limits_{\nu  = e,i} {} \left( {{{\left\| {\partial _t^k{u^\nu }} \right\|}^2} + {{\left\| {\partial _t^{k + 1}{u^\nu }} \right\|}^2} + C{{\left\| {\partial _t^k{\Theta ^\nu }} \right\|}^2}} \right) + C\sum\limits_{\nu  = e,i} {} |||{V^\nu }|||_s^2|||Z||{|_s}. \\
  \end{split}\end{equation*}
This proves \eqref{3.37}.

Next, from 
the first equation of system \eqref{2.3}, we get
$$\partial _t^s{Q^\nu } = \partial _t^s\left( {\frac{1}{{2{\theta
^\nu }}}{{\left| {{u^\nu }} \right|}^2} - {u^\nu } \cdot \nabla
{Q^\nu } - \frac{1}{{{\theta ^\nu }}}{\Theta ^\nu }} \right) -
2\nabla  \cdot
\partial _t^s{u^\nu } - \partial _t^s{u^\nu } \cdot \nabla {\bar
q^\nu }.
$$
Thus, \eqref{3.38} follows from the Leibniz formula and Lemma
\ref{L2.2}. \hfill $\Box$

Using Lemmas \ref{L3.3}-\ref{L3.5}, we get the following result.
\begin{lemma} (Relation of recurrence)
\label{L3.7} 
 Suppose that the conditions of Theorem \ref{T1} hold and $\omega_T$
is small enough independent of $T$, then there exists a positive
constant $C_0$ such that, for all $k \in \N$ and $\alpha \in \N^3$
with $|\alpha| \geq 1$ and $k+|\alpha| \leq s$, we have
\begin{equation}
\label{3.42}
\begin{split}
& \frac{d}{{dt}}\left( {\sum\limits_{\nu  = e,i} {} \sum\limits_{\beta  \le \alpha } {} \left\langle {A_0^\nu({p^\nu },{\theta ^\nu }) V_{k,\beta }^\nu ,V_{k,\beta }^\nu } \right\rangle  + ||{F_{k,\beta }}|{|^2} + ||{G_{k,\beta }}|{|^2}} \right) + {C_0}\sum\limits_{\nu  = e,i} {} ||\partial _t^k{V^\nu }||_{|\alpha |}^2 \\
  \le & C\sum\limits_{\nu  = e,i} {} \left( {||\partial _t^k{V^\nu }||_{|\alpha | - 1}^2 + ||\partial _t^{k + 1}{u^\nu }||_{|\alpha | - 1}^2} \right) + C\sum\limits_{\nu  = e,i} {} |||{V^\nu }|||_s^2|||Z||{|_s}. \\
\end{split}
\end{equation}
\end{lemma}
\noindent{\textbf{Proof.}} For all $k\in \mathbb{N}$, $\beta \in
\mathbb{N}^3$ with $|\beta|\ge 1$ and $k+|\beta| \le s$, it follows
from Lemma \ref{L3.3} and Lemma \ref{L3.5}  that
\begin{equation*}\begin{split}
&\frac{d}{{dt}}\left( {\sum\limits_{\nu  = e,i} {} \left\langle {A_0^\nu ({p^\nu },{\theta ^\nu })V_{k,\beta }^\nu ,V_{k,\beta }^\nu } \right\rangle  + ||{F_{k,\beta }}|{|^2} + ||{G_{k,\beta }}|{|^2}} \right) + {C_0}\sum\limits_{\nu  = e,i} {} ||\left( {u_{k,\beta }^\nu ,Q_{k,\beta }^\nu ,\Theta _{k,\beta }^\nu } \right)|{|^2} \\
  \le & C\sum\limits_{\nu  = e,i} {} \left( {||\partial _t^k\left( {{Q^\nu },{u^\nu },{\Theta ^\nu }} \right)||_{|\beta | - 1}^2 + ||\partial _t^{k + 1}{u^\nu }||_{|\beta | - 1}^2} \right) + C\sum\limits_{\nu  = e,i} {} |||{V^\nu }|||_s^2|||Z||{|_s}, \\
 \end{split}\end{equation*}
The summation of these inequalities for all $\beta$ up to $|\beta|
\le |\alpha|$, together with Lemma \ref{L3.4}, \eqref{3.42} follows.
\hfill $\Box$

\vspace{3mm}

Next, using Lemma \ref{L3.4} and Lemma \ref{L3.6}, we obtain the
following result.
\begin{lemma}
\label{L3.8} 
 Suppose that the conditions of Theorem \ref{T1} hold and $\omega_T$
is small enough independent of $T$, then there exists a positive
constant $C_0$ such that, for all $s \ge 3$, we have
\begin{equation}
\label{3.43}
\begin{split}
&
 \frac{d}{{dt}}\left( {\sum\limits_{\nu  = e,i} {} \left\langle {A_0^\nu ({p^\nu },{\theta ^\nu })\partial _t^s{V^\nu },\partial _t^s{V^\nu }} \right\rangle  + ||\partial _t^sF|{|^2} + ||\partial _t^sG|{|^2}} \right) + {C_0}\sum\limits_{\nu  = e,i} {} ||\partial _t^s{V^\nu }|{|^2} \\
  \le & C\sum\limits_{\nu  = e,i} {} ||\partial _t^{s - 1}{V^\nu }||_{|1}^2 + C\sum\limits_{\nu  = e,i} {} |||{V^\nu
  }|||_s^2|||Z||{|_s}.
\end{split}
\end{equation}
\end{lemma}
%
%


\section{Proof of Theorem \ref{T1}}
The proof of Theorem \ref{T1} is mainly based on the following a
priori estimates which is a consequence on the estimates obtained in
Section 3.
\begin{proposition} ( A
priori estimates.)
\label{P3.1}
  Assume that the conditions of Theorem \ref{T1} holds and $\omega_T$
is sufficiently small independent of $T$, then for all $t \in [0,
T]$, there exist positive constants $C_0$ and $C$ such that
\begin{equation}
\label{3.2}
\begin{split}
|||Z\left( {t, \cdot } \right)|||_s^2 + {C_0}\int_0^t  \big( &
\sum\limits_{e,i}  |||{V^\nu }\left( {\tau , \cdot }
\right)|||_s^2 + |||F\left( {\tau , \cdot } \right)|||_{s - 1}^2 \\
+ & |||\nabla G\left( {\tau , \cdot } \right)|||_{s - 2}^2 +
|||{\partial _t}G\left( {\tau , \cdot } \right)|||_{s -
2}^2\big)d\tau \le C||{Z_0}||_s^2.
\end{split}
\end{equation}
\end{proposition}
\noindent{\textbf{Proof.}} First,  for any fixed index $k\in
\mathbb{N}$ with $k \le s-1$, we employ the induction on space
derivatives $|\alpha|$ with $1 \le |\alpha| \le s-k$ for
\eqref{3.42}. The step of the induction is increasing from $|\alpha|
= 1$ to $|\alpha| = s-k.$ More precisely, for $|\alpha| \ge 2$, $
{|\partial _t^k{V^\nu }||_{|\alpha | - 1}^2} $ on the right-hand
side of \eqref{3.42} can be controlled by $ \sum\limits_{\nu  = e,i}
{} ||\partial _t^k{V^\nu }||_{|\alpha |}^2 $ in the proceeding step
on the left-hand side of \eqref{3.42} multiplying a small positive
constant. Then we have
\begin{equation}
\label{3.44}
\begin{split}
&
 \frac{d}{{dt}}\left( {\sum\limits_{\left| \alpha  \right| \le s - k} {} {m_{k,\alpha }}\left( {\sum\limits_{\nu  = e,i} {} \left\langle {A_0^\nu V_{k,\alpha }^\nu ,V_{k,\alpha }^\nu } \right\rangle  + ||{F_{k,\alpha }}|{|^2} + ||{G_{k,\alpha }}|{|^2}} \right)} \right) + {C_0}\sum\limits_{\nu  = e,i} {} ||\partial _t^k{V^\nu }||_{s - k}^2 \\
  \le & C\sum\limits_{\nu  = e,i} {} \left( {||\partial _t^k{V^\nu }|{|^2} + ||\partial _t^{k + 1}{u^\nu }||_{s - k - 1}^2} \right) + C\sum\limits_{\nu  = e,i} {} |||{V^\nu }|||_s^2|||Z||{|_s}. \\
\end{split}
\end{equation}
where $m_{k,\alpha }$ are positive constants.

Next, we continue to employ  the induction on time derivatives $k$
from $k = s$ to $k = 0$.  The corresponding estimate for $k = s$ is
given by \eqref{3.43}. For $k = s - 1$, \eqref{3.44} yields
\begin{equation}
\label{3.45}
\small\begin{split} &
 \frac{d}{{dt}}\left( {\sum\limits_{\left| \alpha  \right| \le 1} {} {m_{s - 1,\alpha }}\left( {\sum\limits_{\nu  = e,i} {} \left\langle {A_0^\nu V_{s - 1,\alpha }^\nu ,V_{s - 1,\alpha }^\nu } \right\rangle  + ||{F_{s - 1,\alpha }}|{|^2} + ||{G_{s - 1,\alpha }}|{|^2}} \right)} \right) + {C_0}\sum\limits_{\nu  = e,i} {} ||\partial _t^{s - 1}{V^\nu }||_1^2 \\
  \le& C\sum\limits_{\nu  = e,i} {} \left( {||\partial _t^{s - 1}{V^\nu }|{|^2} + ||\partial _t^s{u^\nu }||_{s - k - 1}^2} \right) + C\sum\limits_{\nu  = e,i} {} |||{V^\nu }|||_s^2|||Z||{|_s}. \\
\end{split}
\end{equation}
Obviously, the term $\sum\limits_{\nu  = e,i} {} ||\partial _t^{s -
1}{V^\nu }||_{|1}^2$ on the right-hand side of \eqref{3.43} can be
controlled by the same term on the left-hand side of \eqref{3.45}
multiplying a small positive constant. In a similar way, the term $
\sum\limits_{\nu = e,i} {}||\partial _t^{k + 1}{u^\nu }||_{s - k -
1}^2$ can be controlled by $ \sum\limits_{\nu = e,i} {}  ||\partial
_t^k{V^\nu }|{|^2}$ in the proceeding step. Then by induction on
$k$, we obtain
\begin{equation}
\label{3.46}
\small\begin{split} &
  \frac{d}{{dt}}\left( {\sum\limits_{k + \left| \alpha  \right| \le s} {} {m_{k,\alpha }}\left( {\sum\limits_{\nu  = e,i} {} \left\langle {A_0^\nu V_{k,\alpha }^\nu ,V_{k,\alpha }^\nu } \right\rangle  + ||{F_{k,\alpha }}|{|^2} + ||{G_{k,\alpha }}|{|^2}} \right)} \right) + {C_0}\sum\limits_{\nu  = e,i} {} \sum\limits_{k = 0}^s {} ||\partial _t^k{V^\nu }||_{s - k}^2 \\
  \le& C\sum\limits_{\nu  = e,i} {} \sum\limits_{k = 0}^{s - 1} {} \left( {||\partial _t^k{V^\nu }|{|^2} + ||\partial _t^{k + 1}{u^\nu }|{|^2}} \right) + C\sum\limits_{\nu  = e,i} {} |||{V^\nu }|||_s^2|||Z||{|_s}. \\
\end{split}
\end{equation}
where the positive constants $m_{k, \alpha}$ are possibly amended
based on \eqref{3.44}. Due to the following equivalence relation
$$\sum\limits_{k = 0}^s {} ||\partial _t^k{V^\nu }||_{s - k}^2 \sim
|||V^\nu|||_s^2,$$
then from \eqref{3.25}, \eqref{3.37} and \eqref{3.46}, with a
modification again the constants $m_{k, \alpha}$, we have
\begin{equation*}\begin{split}
& \frac{d}{{dt}}\left( {\sum\limits_{k + \left| \alpha  \right| \le
s} {} {m_{k,\alpha }}\left( {\sum\limits_{\nu  = e,i} {}
\left\langle {A_0^\nu V_{k,\alpha }^\nu ,V_{k,\alpha }^\nu }
\right\rangle  + ||{F_{k,\alpha }}|{|^2} + ||{G_{k,\alpha }}|{|^2}}
\right)} \right)
+ \frac{3}{2}\sum\limits_{\nu  = e,i} {} |||{V^\nu }|||_s^2\\
 \le &
C\sum\limits_{\nu  = e,i} {} |||{V^\nu }|||_s^2|||Z||{|_s}.
\end{split}\end{equation*}
Since $\omega_T$ is sufficiently small, we further get 
\begin{equation*}\begin{split}
\frac{d}{{dt}}\left( {\sum\limits_{k + \left| \alpha  \right| \le s}
{} m_{k,\alpha }\left( {{}\sum\limits_{\nu  = e,i} {} \left\langle
{A_0^\nu V_{k,\alpha }^\nu ,V_{k,\alpha }^\nu } \right\rangle  +
||{F_{k,\alpha }}|{|^2} + ||{G_{k,\alpha }}|{|^2}} \right)} \right)
+ \sum\limits_{\nu  = e,i} {} |||{V^\nu }|||_s^2 \le 0.
\end{split}\end{equation*}
Noting the following equivalence relation
\begin{equation*}\begin{split}
|||Z|||_s^2 \sim {\sum\limits_{k + \left| \alpha  \right| \le s} {}
m_{k,\alpha }\left( {\sum\limits_{\nu  = e,i} {} {}\left\langle
{A_0^\nu V_{k,\alpha }^\nu ,V_{k,\alpha }^\nu } \right\rangle  +
||{F_{k,\alpha }}|{|^2} + ||{G_{k,\alpha }}|{|^2}} \right)},
\end{split}\end{equation*}
we obtain
\begin{equation}
\label{3.47}
\small\begin{split} &
  |||Z(t, \cdot )|||_s^2 +  \sum\limits_{\nu  = e,i} {}\int_0^t {}
|||{V^\nu }(\tau , \cdot )|||_s^2d\tau  \le |||Z(0, \cdot )|||_s^2
\le C||{Z_0}||_s^2, \quad t \in [0,T],
\end{split}
\end{equation}
where Lemma \ref{L2.3} is used.

From the second and the last two equations in \eqref{2.3}, we obtain
the estimates for electric-magnetic fields as
\begin{equation}
\label{3.48}
\small\begin{split} &
  |||F|||_{s - 1}^2 \le C\sum\limits_{\nu  = e,i} {} |||{V^\nu
}|||_s^2 + C\sum\limits_{\nu  = e,i} {} |||{V^\nu
}|||_s^2|||Z||{|_s},
\end{split}
\end{equation}
and
\begin{equation}
\label{3.49}
\small\begin{split} &
  |||{\partial _t}G|||_{s - 2}^2 + |||{\nabla _x}G|||_{s - 2}^2 \le
C\sum\limits_{\nu  = e,i} {} |||{V^\nu }|||_s^2 + C\sum\limits_{\nu
= e,i} {} |||{V^\nu }|||_s^2|||Z||{|_s}.
\end{split}
\end{equation}
Due to the fact that $\omega_T$ is sufficiently small, combining
\eqref{3.47}-\eqref{3.49} yields \eqref{3.2}.
 \hfill $\Box$

\vspace{3mm}

 It is obvious that \eqref{3.2} gives \eqref{1.12} and
the global existence of smooth solution $(n^\nu, u^\nu, \theta^\nu,
E, B)$ to periodic problem \eqref{1.1*}-\eqref{1.3}. Finally, for
all $k\in \N$ and $\beta \in \N^3$ with $k + |\beta| \le s-1,$ from
\eqref{1.12}, we get
$$\partial _t^k{\partial ^\beta }\left( {{n^\nu } - {{\bar n}^\nu
},{u^\nu },{\theta ^\nu } - 1,E - \bar E} \right) \in {L^2}\left(
{{\R^ + };{L^2}(\T)} \right) \cap {W^{1,\infty }}\left( {{\R^ +
};{L^2}(\T)} \right), \quad \nu = e, i,$$
which implies \eqref{1.13}-\eqref{1.14}. Furthermore, if $k +
|\beta| \ge 1,$ noticing $\bar B$ is a constant vector, we have
$$
\partial _t^k{\partial ^\beta }B \in {L^2}\left( {{\R^ + };{L^2}(\T)}
\right) \cap {W^{1,\infty }}\left( {{\R^ + };{L^2}(\T)} \right),$$
which gives \eqref{1.15}.

 \hfill $\Box$


\vspace{3mm}

 \noindent {\sc Acknowledgments~:}
The authors are grateful to the referee for the comments. The first
author would like to  express his sincere gratitude to Professor
Yue-Jun Peng of Universit\'e Blaise Pascal for excellent directions
in France. The authors are supported by the the BNSF (1164010,
1132006), NSFC (11371042), the key fund of the Beijing education
committee of China,  the general project of scientific research
project of the Beijing education committee of China, NSF of Qinghai
Province, the collaborative innovation center on Beijing
society-building and social governance, the China postdoctoral
science foundation funded project, the Project supported by Beijing
Postdoctoral Research Foundation, the government of Chaoyang
district postdoctoral research foundation, the 2016 Beijing project
of scientific activities for the excellent students studying abroad
and the Beijing University of Technology foundation funded project.



\vspace{4mm}


\begin{thebibliography}{99}




\bibitem{Ch84} F.~Chen,
Introduction to Plasma Physics and Controlled Fusion. Vol. 1, {\it
Plenum Press, New York.} 1984.

\bibitem{CJW00} G.Q.~Chen, J.W.~Jerome, D.H.~Wang,
Compressible Euler-Maxwell equations, {\it Transport Theory and
Statistical Physics}, 29 (2000) 311-331.


\bibitem{DDS12} P.~Degond, F.~Deluzet, D.~Savelief,
Numerical approximation of the Euler-Maxwell model in the
quasineutral limit, {\it J.Comput. Phys.} 231 (2012), 1917-1946.

\bibitem{Du11} R.J.~Duan,
Global smooth flows for the compressible Euler-Maxwell system~: the
relaxation case, {\it J.Hyper. Diff. Equations}, 8 (2011) 375-413.

\bibitem{DLZ12} R.J.~Duan, Q.~Q.~Liu, C.J.~Zhu,
The Cauchy problem on the compressible two-fluids Euler-Maxwell
equations, {\it SIAM J. Math. Anal.} 44 (2012) 102-133.

\bibitem{Ev98} L.C.~Evans,
Partial Differential Equations, Graduate Studies in Mathematics, 19,
American Mathematical Society, Providence, RI, 1998.

\bibitem{FPW15} Y.H.~Feng, Y.J. Peng, S.~Wang,
Stability of non-constant equilibrium solutions for two-fluid
Euler-Maxwell systems, {\it Nonlinear Anal. Real World}, 26 (2015)
372-390.

\bibitem{FWK11} Y.H.~Feng, S.~Wang, S.~Kawashima,
Global Existence and Asymptotic Decay of Solutions to the
Non-Isentropic Euler-Maxwell System, {\it Math. Models Methods Appl.
Sci.} 24 (2014)
 2851-2884.


\bibitem{FWL16a} Y.H.~Feng, S.~Wang, X.~Li,
Stability of non-constant steady-state solutions for non-isentropic
Euler-Maxwell system  with a temperature damping term, {\it Math.
Methods  Appl. Sci. } 39 (2016) 2514-2528.




\bibitem{GM11} P.~Germain, N.~Masmoudi,
Global existence for the Euler-Maxwell system, preprint,
arXiv:1107.1595, 2011.


\bibitem{GIP16} Y.~Guo, A.D. Ionescu, B. Pausader,
Global solutions of the Euler-Maxwell two-fluid system in 3D,
preprint, Ann. Of Math. 183 (2016), 377-498.

\bibitem{GS05} Y.~Guo, W.~Strauss,
Stability of semiconductor states with insulating and contact
boundary conditions, {\it Arch. Ration. Mech. Anal.} 179 (2005)
1-30.



\bibitem{Ka75} T.~Kato,
The Cauchy problem for quasi-linear symmetric hyperbolic systems,
{\it Arch. Ration. Mech. Anal.} 58 (1975) 181--205.

\bibitem{KM81} S.~Klainerman, A.~Majda,
Singular limits of quasilinear hyperbolic systems with large
parameters and the incompressible limit of compressible fluids, {\it
Comm. Pure Appl. Math.} 34 (1981) 481-524.

\bibitem{LWF16a} X. Li, S. Wang, Y.H. Feng, Stability of non-constant steady-state solutions for bipolar non-isentropic
 Euler-Maxwell equations with damping terms,  {\it
Z. Angew. Math. Phys.}  67(5) (2016) 1-27.

\bibitem{LP17} C. Liu, Y.J. Peng, Stability of periodic steady-state solutions to a non-isentropic
Euler-Maxwell system. Preprint 2017.

\bibitem{LZ13} Q.Q. Liu, C.J. Zhu, Asymptotic stability of stationary solutions to
the compressible
 Euler-Maxwell equations,  {\it
Indiana. Univ. Math. J.}  62 (2013) 1203-1235.



\bibitem{Ma84} A.~Majda, Compressible Fluid Flow and Systems of Conservation Laws in
Several Space Variables, {\it Springer-Verlag, New York,} 1984.



\bibitem{MRS90} P.A.~Markowich, C.A.~Ringhofer, C.~Schmeiser,
Semiconductor Equations, Springer-Verlag, New York, 1990.


\bibitem{MN79} A.~Matsumura, T.~Nishida,
{ The initial value problem for the equation of motion of
compressible viscous and heat-conductive fluids.}, {\it Proc. Japan
Acad, Ser A,}  55 (1979) 337-342.

\bibitem{MN80} A.~Matsumura, T.~Nishida,
{ The initial value problem for the equation of motion of viscous
and heat-conductive gases.} {\it J. Math. Kyoto Univ}, 20 (1980)
67-104.





\bibitem{PW08a} Y.J.~Peng, S.~Wang,
Convergence of compressible Euler-Maxwell equations to
incompressible Euler equations, {\it Commun. Part. Differ.
Equations,} 33 (2008) 349-376.

\bibitem{PWG11} Y.J.~Peng, S.~Wang, G.~L.~Gu,
Relaxation limit and global existence of smooth solutions of
compressible Euler-Maxwell equations, {\it SIAM J. Math. Anal.} 43
(2011) 944-970.

\bibitem{Pe12} Y.J.~Peng,
Global existence and long-time behavior of smooth solutions of
two-fluid Euler-Maxwell equations, {\it Ann. I. H. Poincare-AN.} 29
 (2012) 737-759.

\bibitem{Pe13} Y.J.~Peng,
Stability of non-constant equilibrium solutions for Euler-Maxwell
equations, {\it J. Math. Pures Appl.} 103 (2015) 39-67.

\bibitem{RG69} H.~Rishbeth, O.K.~Garriott,
Introduction to Ionospheric Physics. {\it Academic Press,} London,
1969.

\bibitem{UK11} Y.~Ueda, S.~Kawashima,
Decay property of regularity-loss type for the Euler-Maxwell system,
{\it Methods Appl. Anal.} 18 (2011) 245-267.

\bibitem{UWK12} Y. Ueda, S. Wang, S.~Kawashima,
Dissipative structure of the regularity-loss type and time
asymptotic decay of solutions for the Euler-Maxwell system, {\it
SIAM J. Math. Anal.} 44 (2012) 2002-2017.

\bibitem{WFL12} S.~Wang, Y.H.~Feng, X.~Li,
The asymptotic behavior of globally smooth solutions of bipolar
non-isentropic compressible Euler-Maxwell system for plasma, {\it
SIAM J. Math. Anal.} 44 (2012) 3429-3457.


\bibitem{WFL14}
 S. Wang, Y.H. Feng, X. Li,  The asymptotic behavior of globally smooth solutions of
non-isentropic  Euler-Maxwell  equations  for  plasmas. {\it Appl.
Math. Comput}, 231 (2014) 299-306.

\bibitem{WZ14} H.Y. Wen, C.J. Zhu,
Global symmetric classical solutions of the full compressible
Navier-Stokes equations with vacuum and large initial data. {\it J.
Math. Pures Appl}, 102 (2014) 498-545.

\bibitem{Xu11} J.~Xu,
Global classical solutions to the compressible Euler-Maxwell
equations, {\it SIAM J. Math. Anal.} 43 (2011) 2688-2718.

\bibitem{Yo04} W.A.~Yong,
Entropy and global existence for hyperbolic balance laws, {\it Arch.
Ration. Mech. Anal.} 172 (2004) 247-266.
\end{thebibliography}
\end{document}